\pdfoutput=1
\RequirePackage{ifpdf}
\ifpdf 
\documentclass[pdftex]{sigma}
\else
\documentclass{sigma}
\fi

\numberwithin{equation}{section}

\newtheorem{Theorem}{Theorem}[section]

\newtheorem{Lemma}[Theorem]{Lemma}
\newtheorem{prop}[Theorem]{Proposition}
\newtheorem{Conjecture}[Theorem]{Conjecture}
{\theoremstyle{definition}

\newtheorem{Remark}[Theorem]{Remark}
}

\newcommand{\CP}{\mathbb{C}P}

\newcommand{\R}{\mathbb{R}}
\newcommand{\RR}{\mathbb{R}}
\newcommand{\C}{\mathbb{C}}
\newcommand{\CC}{\mathbb{C}}
\newcommand{\Z}{\mathbb{Z}}
\newcommand{\ZZ}{\mathbb{Z}}

\newcommand{\CL}{\mathcal{L}}
\newcommand{\CM}{\mathcal{M}}

\newcommand{\Hom}{{\rm Hom}}

\begin{document}

\allowdisplaybreaks

\renewcommand{\thefootnote}{$\star$}

\renewcommand{\PaperNumber}{053}

\FirstPageHeading

\ShortArticleName{Examples of Matrix Factorizations from SYZ}

\ArticleName{Examples of Matrix Factorizations from SYZ\footnote{This
paper is a contribution to the Special Issue ``Mirror Symmetry and Related Topics''. The full collection is available at \href{http://www.emis.de/journals/SIGMA/mirror_symmetry.html}{http://www.emis.de/journals/SIGMA/mirror\_{}symmetry.html}}}

\Author{Cheol-Hyun CHO, Hansol HONG and Sangwook LEE}

\AuthorNameForHeading{C.-H.~Cho, H.~Hong and S.~Lee}

\Address{Department of Mathematics, Research Institute of Mathematics, Seoul National University,\\
 1 Kwanak-ro, Kwanak-gu, Seoul, South Korea}
\Email{\href{mailto:chocheol@snu.ac.kr}{chocheol@snu.ac.kr}, \href{mailto:hansol84@snu.ac.kr}{hansol84@snu.ac.kr}, \href{mailto:leemky7@snu.ac.kr}{leemky7@snu.ac.kr}}

\ArticleDates{Received May 15, 2012, in f\/inal form August 12, 2012; Published online August 16, 2012}

\Abstract{We f\/ind matrix factorization corresponding to an anti-diagonal in $\CP^1 \times \CP^1$, and circle f\/ibers in weighted projective lines using
the idea of Chan and Leung of Strominger--Yau--Zaslow transformations. For the tear drop orbifolds,
we apply this idea to f\/ind  matrix factorizations for two types of potential, the usual Hori--Vafa potential or
the bulk deformed (orbi)-potential.
We also show that the direct sum of anti-diagonal with its shift, is equivalent to the direct sum of
central torus f\/ibers with holonomy $(1,-1)$ and $(-1,1)$ in the Fukaya category of $\CP^1 \times \CP^1$, which was predicted by
Kapustin and Li from B-model calculations.}


\Classification{53D37; 53D40; 57R18}

\renewcommand{\thefootnote}{\arabic{footnote}}
\setcounter{footnote}{0}

\section{Introduction}

The Strominger--Yau--Zaslow (SYZ for short) \cite{SYZ} conjecture provides a geometric way to understand mirror symmetry phenomenons. Recently, Chan and Leung~\cite{CL} have shown that in~$\CP^1$
the Lagrangian Floer chain complex  between the equator and a generic Lagrangian torus f\/iber   corresponds by SYZ to the matrix factorization of the Landau--Ginzburg (LG for short)  superpotential~$W$.  The general idea is as follows.  To f\/ind a matrix factorization corresponding to a Lagrangian submanifold, say $L_0$, they
consider a family of Floer chain complex of the pair $(L_u, L_0)$ for all torus f\/ibers~$L_u$
(with all possible holonomies), and use the information of holomorphic strips (which
varies as $L_u$ changes) and apply SYZ transformation to construct the matrix factorization for~$L_0$
(see Section~\ref{section2} for more details).

Their observation is very insightful to understand the homological mirror symmetry between Lagrangian Floer theory of toric manifolds and matrix factorization of LG superpotential~$W$, but the procedure is known to work only for a $\CP^1$ (or a product of $\CP^1$ with Lagrangian submanifold given by
product of equators).  They also found the corresponding matrix factorization for $\CP^2$, but the description of Floer chain complex is not complete.

 In this paper, we provide more evidence on this correspondence following their ideas. The f\/irst
 new example is the case of the anti-diagonal Lagrangian submanifold  in the symplectic manifold $\CP^1 \times \CP^1$. In fact, Kapustin and Li already conjectured  in \cite{KL} that the anti-diagonal should correspond to  a specif\/ic matrix factorization of LG superpotential $W = x + \frac{q}{x} + y + \frac{q}{y}$, and we verify this conjecture using this procedure.

\begin{prop}\label{prop1}
For $\CP^1 \times \CP^1$, the anti-diagonal Lagrangian submanifold corresponds to the following matrix factorization by SYZ transformation $($in the sense of {\rm \cite{CL})}
 \begin{gather*}
(x+y)\left(1 + \dfrac{q}{xy} \right) =   x + \frac{q}{x} + y + \frac{q}{y}.
\end{gather*}
\end{prop}

  For this,  we deform  generic Lagrangian torus f\/ibers into  specif\/ic forms via Hamiltonian isotopy, and analyze the
Floer cohomology of the anti-diagonal, with  ``deformed'' generic torus f\/ibers, and apply Chan and Leung's SYZ transformation to
f\/ind the corresponding matrix factorization.

From $B$-model calculations, Kapustin and Li further conjectured in~\cite{KL} that in the Fukaya category, the direct sum of
anti-diagonal $A$ and its shift $A[1]$, is isomorphic to the direct sum of two f\/ibers $T_{1,-1}$ and $T_{-1,1}$ of holonomy
$(1,-1)$ and $(-1,1)$ respectively.
We also verify this conjecture by computing the Floer cohomology and  products  between these objects, and f\/inding
a homomorphism which induces this isomorphism.

\begin{Theorem}[Theorem~\ref{thm:equi}]
In the derived Fukaya category of $\CP^1 \times \CP^1$, $A \oplus A[1]$ is  equivalent  to $T_{1,-1} \oplus T_{-1,1}$.
\end{Theorem}

Namely, $\CP^1 \times \CP^1$
has four Lagrangian torus f\/ibers, whose Floer cohomology groups are non-vanishing.
It is given by the central f\/iber $T^2$, with holonomy $(1,1)$, $(1,-1)$, $(-1,1)$, $(-1,-1)$.
Central f\/iber with holonomy $(1,-1)$ (or $(-1,1)$), which we denote by $T_{1,-1}$ (or $T_{-1,1}$) has va\-ni\-shing $m_0$, and
hence is unobstructed. The anti-diagonal Lagrangian submanifold $A$, is monotone Lagrangian submanifold of minimal Maslov index~4, hence unobstructed. Hence the
Lagrangian Floer cohomology among $\{T_{1,-1}, T_{-1,1}, A, A[1]\}$ can be def\/ined,
where $A[1]$ is regarded as an object of (derived) Fukaya category. In Section~\ref{section5}, we compute the Floer cohomology
between these objects as well as several~$m_2$ products between them to verify the conjecture.

\looseness=-1
Our second type of examples are  weighted projective lines, which are toric orbifolds. There is an interesting new phenomenon due to {\em bulk deformation by twisted sectors of toric orbifolds}.
First of all, these weighted~$\CP^1$'s have  Landau--Ginzburg mirror superpotential $W:\C^* \to \C$, and
the f\/irst author and Poddar has recently
developed a Lagrangian Floer cohomology theory for toric orbifolds, and superpotential $W$ can be def\/ined from
the data of smooth holomorphic discs in toric orbifolds. We consider the Floer chain complex of a central f\/iber of the weighted~$\CP^1$ and a generic torus f\/iber, and from this we can f\/ind the corresponding matrix factorization of~$W$.

\begin{prop}
For a weighted $\CP^1$ with $\ZZ/ m\ZZ$-singularity on the left and $\ZZ/ n\ZZ$ on the right, the central fiber corresponds to the following matrix factorization by SYZ transformation:
\begin{gather*}
 \left( 1-\dfrac{z}{\alpha q} \right) \left( \displaystyle\sum_{k=0}^{n} \frac{q^{ \frac{m}{n}  k}}{\alpha^{k}} \left(\frac{q^{\frac{m+n}{n}}}{z} \right)^{n-k} -\sum_{k=1}^{m}\alpha^k q^{k}z^{m-k} \right)
 = z^m+\dfrac{q^{m+n}}{z^n}-\left( \alpha^m q^{m}+\dfrac{q^{m}}{\alpha^n} \right).
\end{gather*}
\end{prop}

(See Sections \ref{SecTO} and \ref{SecWPR}.)

Then, we can turn on bulk deformation $\frak{b}$ by twisted sectors to obtain a bulk deformed
mirror superpotential $W^{\frak{b}}$. This potential has additional terms from the data
of orbifold holomorphic discs in toric orbifolds. Once bulk deformation~$\frak{b}$ is chosen
(so that there is a torus f\/iber~$L$ whose Floer cohomology is non-vanishing), then
we consider the Floer chain complex of $L$ with a generic torus f\/iber to f\/ind
a matrix factorization of bulk deformed LG superpotential~$W^{\frak{b}}$.
In this case, we f\/ind the corresponding matrix factorization
of~$W^{\frak{b}}$ by additionally considering orbifold holomorphic strips.

\section{Preliminaries}\label{section2}

We recall Strominger--Yau--Zaslow conjecture brief\/ly. The classical form of mirror symmetry considers mirror pairs of Calabi--Yau 3-folds $X$ and $\check{X}$, and the symplectic geometry (Gromov--Witten invariants) of $X$ corresponds to  complex geometry (periods) of $\check{X}$.  The SYZ conjecture is, roughly speaking, a geometric tool to f\/ind the mirror manifold, as a dual torus f\/ibration. We state the conjecture in the following form from~\cite{Gross} (see also~\cite{Au}):

\begin{Conjecture}[\cite{SYZ}] If two Calabi--Yau $n$-folds $X$ and $\check{X}$ are mirror to each other, then there exist special Lagrangian fibrations $f:X\to B$ and $\check{f}:\check{X}\to B$, whose generic fibers are tori. Furthermore, these fibrations are dual, namely $X_b=H^1(\check{X}_b,\RR/\ZZ)$ and
\[
\check{X}_b=H^1(X_b,\RR/\ZZ),
\] when $X_b$ and  $\check{X}_b$ are nonsingular torus fibers over $b\in B$.
\end{Conjecture}

Toric Fano manifolds $X$, which are torus f\/ibrations over the moment polytopes, has a mirror given by
a Landau--Ginzburg model  $(\check{X},W)$. Torus f\/ibers become singular over the facets of the moment polytope, and the singularity of the f\/ibration is measured by the Landau--Ginzburg superpotential~$W$, which can be constructed from the Maslov index two holomorphic discs in~$X$ with boundary on torus
f\/ibers~\cite{CO,FOOO}.  The homological mirror symmetry due to Kontsevich (in this setting) asserts that
the derived Fukaya category ${\rm DFuk}(X)$ of a toric Fano manifold $X$ is equivalent, as a triangulated category, to the category of matrix factorizations $MF(\check{X},W)$ of the mirror
 Landau--Ginzburg model  $(\check{X},W)$. The latter category is equivalent to the
 category of singularites $D_{\rm Sg}(\check{X},W)$ (see Orlov~\cite{Or}).

A matrix factorization of a Landau--Ginzburg model $(\check{X},W)$ is a square matrix $M$ of
even dimensions with entries in the coordinate ring $\C[\check{X}]$ and of the form
\begin{gather*}
M =  \begin{pmatrix} 0 & F \\ G & 0 \end{pmatrix},
\end{gather*}
such that
\begin{gather*}
M^2 = (W -\lambda){\rm Id}
\end{gather*}
for some $\lambda \in \C$. It is well-known that $M$ is a non-trivial element of $MF(\check{X},W)$ only
if  $\lambda$ is a~critical value of $W$ (see~\cite{Or}).

The idea of Chan and Leung will be explained in more detail in the next section, but we f\/irst explain
the Lagrangian Floer theory behind this correspondence.
Let $L_0$, $L_1$ be a Lagrangian submanifold in a symplectic manifold $(X,\omega)$.
Let $J$ be a compatible almost complex structure. One considers
$J$-holomorphic discs $u:(D^2,\partial D^2) \to (X,L)$ with Lagrangian boundary conditions,
and denote by $\CM_k(L,\beta)$ be the moduli space of such $J$-holomorphic discs of homotopy
class $\beta \in \pi_2(X,L)$ with $k$ boundary marked points. We denote by~$\mu(\beta)$
the Maslov index of~$\beta$. The dimension of the moduli space is given by $n + \mu(\beta) + k -3$.

We further assume that $L_0$, $L_1$ are {\em positive} in the sense that
any non-constant $J$-holomorphic discs have positive Maslov index.
In particular, this implies that the (virtual) dimension of $\CM_1(L,\beta)$ is always at least~$n$
if $\beta \neq 0$ and non-empty. And $\dim (\CM_1(L,\beta)) =n$ exactly when $\mu(\beta)=2$.

We def\/ine the Novikov ring
\[
\Lambda = \left\{ \left. \sum a_i T^{\lambda_i} \right| a_i \in \C, \, \lambda_i \in \R, \, \lim_i \lambda_i = \infty \right\}.
\]
The Lagrangian Floer chain complex $CF(L_0,L_1)$ is generated by intersection points $L_0 \cap L_1$
with coef\/f\/icients $\Lambda$, and its dif\/ferential is def\/ined by
\[
 m_1(\langle p \rangle ) = \sum_{q} n_\alpha(p,q) \langle q \rangle  T^{\omega(\alpha)},
 \]
where the sum is over all $q \in L_0 \cap L_1$, and $n_\alpha(p,q)$ is the
count of isolated $J$-holomorphic strips with boundary on~$L_0$,~$L_1$ (modulo translation action)
of homotopy class $\alpha$, and $\omega(\alpha)$ is the area of such $J$-holomorphic strip.
Such isolated strips have Maslov--Viterbo index one.
We refer readers to \cite{FOOO, Oh} for details.

In general, $m_1^2 \neq 0$ and hence the Lagrangian Floer cohomology cannot be def\/ined in general.
With the above positivity assumption, we have the following {\em Floer complex equation}
\begin{gather*}
 m_1^2 = (W_{L_1} - W_{L_0}) {\rm Id},
 \end{gather*}
where $W_{L_i}$ is def\/ined as follows:
From the evaluation map ${\rm ev}_{0,\beta}:\CM_1(L,\beta) \to L$ at the marked point,
if $\mu(\beta)=2$, the image of ${\rm ev}_{0,\beta}$ is a multiple of fundamental class $[L]$
\[
{\rm ev}_{0,\beta}(\CM_1(L,\beta)) = c_\beta [L]
\]
as it is of dimension $n$, and $\beta$ is of minimal Maslov index.
Then we def\/ine
\begin{gather}\label{W}
W_{L} := \sum_{\mu(\beta)=2} c_\beta T^{\omega(\beta)}.
\end{gather}
\begin{figure}[t]
\centering
\includegraphics[height=2in]{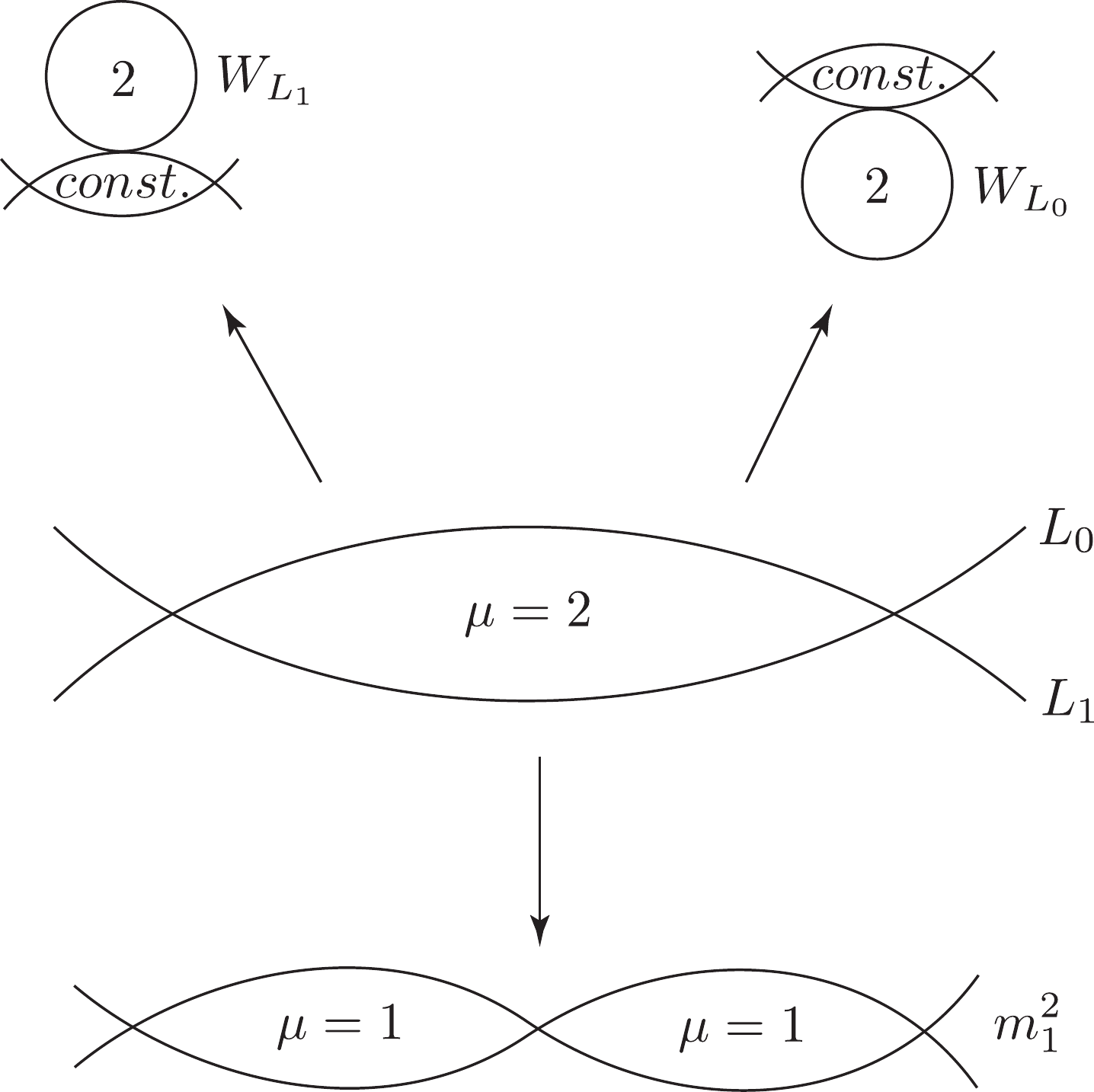}
\caption{Degeneration of index 2 strips.}\label{ind2}
\end{figure}
The Floer complex equation  $m_1^2 = (W_{L_1} - W_{L_0}){\rm Id}$, is obtained by analyzing the moduli space of holomorphic strips of
Maslov--Viterbo index two (see Fig.~\ref{ind2}). Some sequences of $J$-holomorphic strips of Maslov--Viterbo index two, can
degenerate into broken $J$-holomorphic strips, each of which has index one, which contributes to~$m_1^2$. Some sequence of $J$-holomorphic strips of index two can also degenerate into a constant strip
together with a bubble holomorphic disc attached to either upper or lower boundary of
the strip. Discs attached to upper (resp.\ lower) boundary contributes to~$W_{L_1}$ (resp.~$W_{L_0}$)
and it gives~${\rm Id}$ map since the $J$-holomorphic strip is constant.

In fact, one needs to  considers Lagrangian submanifolds $L_0$, $L_1$, equipped with f\/lat line bundles $\CL_0 \to L_0$, $\CL_1 \to L_1$, and the above setting can be extended to this setting. In such a~case,
we put an additional contribution of holonomy
${\rm hol}_{\CL_i}(\partial \beta)$ for each $\beta$ in~\eqref{W}.

Chan and Leung's idea is to compare the Floer complex equation and
that of matrix facto\-ri\-zation $ M^2 = (W -\lambda){\rm Id}$ (via their Fourier transform).
For this, we take~$L_0$ to be a f\/ixed torus f\/iber (corresponding to the critical value $\lambda$)
and vary~$L_1$ as generic torus f\/ibers with holonomy to obtain $W$ as a function on
the mirror manifold.

\section[Chan-Leung's construction for $\CP^1$]{Chan--Leung's construction for $\boldsymbol{\CP^1}$}\label{chanleungcp1}
We recall the result of Chan--Leung~\cite{CL}  in the case of $X=\CC P^1$ for readers' convenience. Recall that $\check{X}=\CC^*$, and the Landau--Ginzburg superpotential is $W=z+\frac{q}{z}$ where $q=T^{t}$ when $[0,t]$ is the moment polytope of $X$. ($W$ can  obtained from the disc potential $e^x T^u + e^{-x} T^{t-u}$ by substituting $z=e^x T^u$).

 By removing north (N) and south (S) pole of $X$,
we regard $X\setminus \{N, S\}$ as a circle f\/ibration over $(0,t)$, and denote by $u$ the
coordinate in $(0,t)$, by $y$ that of the f\/iber circle.  Then the standard symplectic form $\omega$ equals $du \wedge dy$ on  $X\setminus \{N, S\}$.

An equator with trivial holonomy (f\/iber at $t/2$) has non-trivial Floer cohomology, and it
corresponds to the critical point $\sqrt{q}$ of $W$.
By homological mirror symmetry, this should correspond to a skyscraper sheaf at the critical point,
and by Orlov's result, we have a matrix factorization corresponding to it.
The critical value of $W$ is $2\sqrt{q}$ and the corresponding factorization of $W-2\sqrt{q}$ is given in matrix form as
\begin{gather*}
 \begin{pmatrix}
 0 & z-\sqrt{q} \\ 1-\frac{\sqrt{q}}{z} & 0
 \end{pmatrix}.
 \end{gather*}

Chan--Leung's idea is to recover this matrix factorization from the geometry of torus($S^1$) f\/ibration. Let $L_0$ be the central f\/iber with trivial holonomy. We deform $L_0$ to $\tau:[0,3]\to X$ as follows: (in $(u,y)$ coordinate)
\begin{gather*}
\tau(s)=
\begin{cases} ((1-s)t/2,0) & \text{if $0 \leq s \leq 1$,}
\\
((s-1)t/2,0) &\text{if $1 \leq s \leq 2$,}
\\
(t/2,2\pi(s-2)) & \text{if $2 \leq s \leq 3$.}
\end{cases}
\end{gather*}

\begin{figure}[t]
\centering
\includegraphics[height=1.3in]{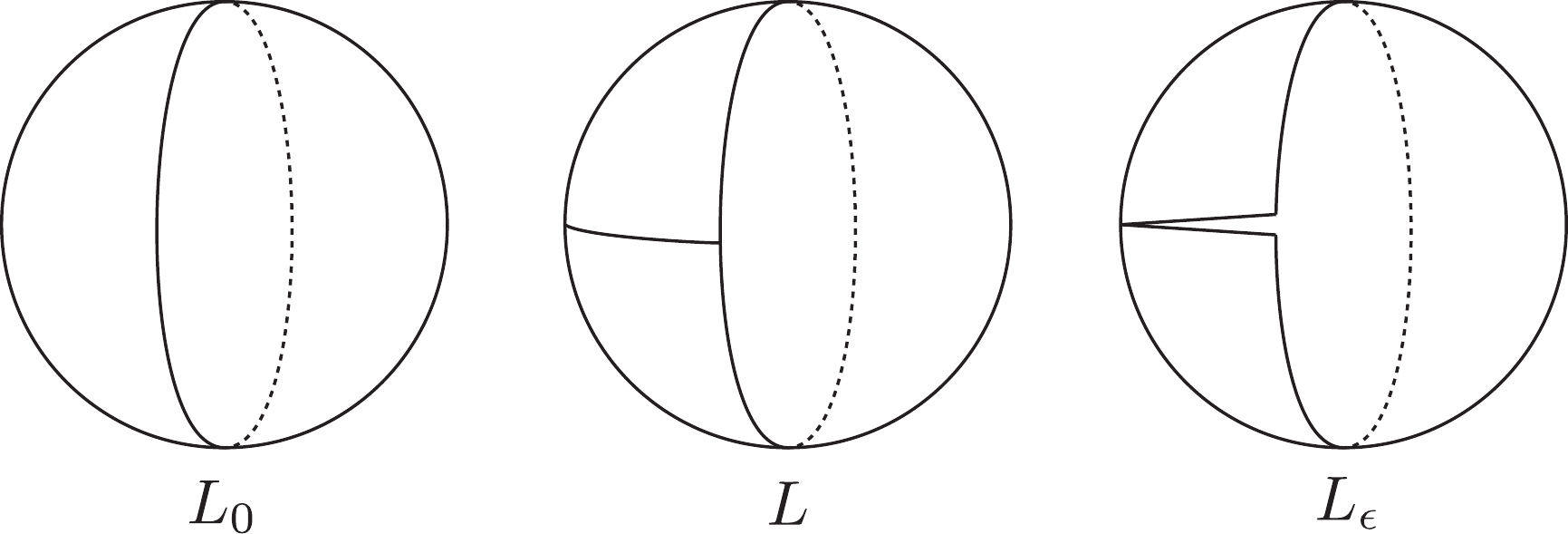}
\caption{Hamiltonian defomation of $L$ and the spike.}\label{spike}
\end{figure}

Namely, $\tau$ f\/irst goes along the zero section from center to the left pole, comes back to the center and at last turns around $L_0$. Let this deformation be denoted by $L$. Note that $L$ still splits $X$ into two equal halves. Since $L$ is too singular, we slightly perturb $L$ to $L_{\epsilon}$ so that it is
still area bisecting, and hence Hamiltonian isotopic to central f\/iber. It is helpful to think of $L$ as a limit of $L_{\epsilon}$ (see Fig.~\ref{spike}).

For each $u \in (0,t)$, we have a corresponding f\/iber $L_u$. Then $L_u$ and $L_{\epsilon}$ meet at two points~$a$ and~$b$, and there occur four holomorphic strips between them. Let $[w]$ be a  homotopy class of a holomorphic strip $w$ between~$a$ and~$b$, namely $[w] \in \pi_2(X;L,L_u;a,b)$. Let $\partial_-[w]$ be the boundary of $w$ on $L_u$. Taking the limit of $L_{\epsilon}$, $\partial_-[w]$ is identif\/ied as an element of $\pi_1(L_u)$.

\begin{figure}[t]
\centering
\includegraphics[height=3.3in]{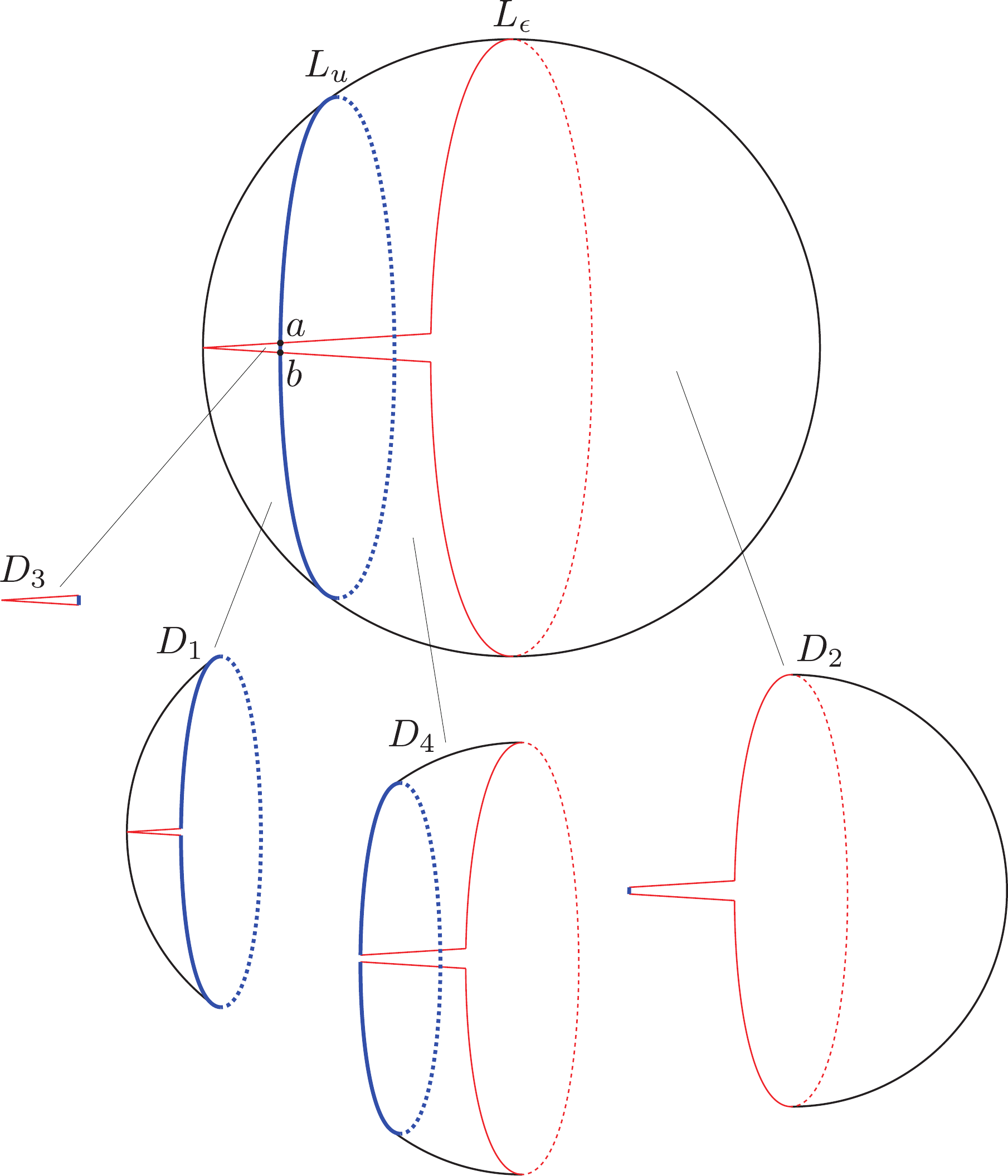}
\caption{Disk splittings in $\CC P^1$.}
\label{dsccp1}
\end{figure}

Now we def\/ine a function $\Psi^{a,b}_L:(0,t/2)\times \pi_1(L_u) \to \RR$ as follows:
\begin{gather}\label{eq:matfun}
\Psi^{a,b}_L(u,[\gamma])=\sum_{\stackrel{[w]\in \pi_2(X;L,L_u;a,b)}{\partial_-(w)=[\gamma]}}\pm n([w]) \exp( \textrm{area($w$)} )  {\rm hol} ([\gamma]) ,
\end{gather}
where the sign is due to the orientation of $w$, and $n([w])$ is the number of holomorphic discs representing $[w]$.

Note that the area of a Maslov index 2 holomorphic disc whose boundary is a toric f\/iber $L_u$ is just given as $u$ or $(t-u)$ (up to times $2\pi$). If we identify $\ZZ \simeq \pi_1(L_u)$, then we have a complete list of (\ref{eq:matfun}):
\begin{gather*}
\Psi^{a,b}_L(u,v)=
\begin{cases} \exp(u) & \text{if $v=1$,}
\\
-\exp({t/2})= -\sqrt{q}&\text{if $v=0$,}
\\
0 & \text{otherwise,}
\end{cases}
\\
\Psi^{b,a}_L(u,v)=
\begin{cases} \exp(0)=1 & \text{if $v=0$,}
\\
-\exp({t/2-u}) =-\dfrac{\sqrt{q}}{\exp(u)} &\text{if $v=-1$,}
\\
0 &\text{otherwise.}
\end{cases}
\end{gather*}

The correspondence of the function values and holomorphic strips is given as follows. $D_i$ are drawn in the above (Fig.~\ref{dsccp1})
\begin{gather*}
\Psi^{a,b}_L(u,1)   \longleftrightarrow   D_1, \qquad\!
\Psi^{a,b}_L(u,0)   \longleftrightarrow   D_2, \qquad\!
\Psi^{b,a}_L(u,0)   \longleftrightarrow   D_3, \qquad\!
\Psi^{b,a}_L(u,-1)   \longleftrightarrow   D_4.
\end{gather*}
Observe that the areas of discs are computed in the limit $L$.

With these functions we make a matrix-valued function $\Psi_L$ by
\begin{gather}\label{eq:mat}
\Psi_L(u,v)= \begin{pmatrix} 0 & \Psi^{a,b}_L(u,v) \\ \Psi^{b,a}_L(u,v) & 0 \end{pmatrix}. \end{gather}

Finally, Fourier transform of \eqref{eq:mat} following \cite{CL2} can be obtained. Each entry of \eqref{eq:mat} is a~function of the form $f=f_v \exp({\langle  u,v \rangle })$, and for such a function $f$ we def\/ine Fourier transform of $f$ as
\[
\hat{f}:= \sum_{v \in \ZZ} f_v \exp({\langle u,v \rangle }) {\rm hol} (v) .
\]
Since ${\rm hol} (v) = \exp({i y v})$, if we adopt a complex coordinate $z=\exp({u+iy})$, then
\[
\hat{f} = \sum_{v \in \ZZ} f_v z^v.
\]
After the Fourier transform, we have
\[
\Psi_L (z)= \begin{pmatrix} 0 & z-\sqrt{q} \\ 1-\frac{\sqrt{q}}{z} & 0 \end{pmatrix},
\]
which is the desired factorization of $W-2\sqrt{q}.$

\section[Anti-diagonal $A$  in $\CC P^1 \times \CC P^1$]{Anti-diagonal $\boldsymbol{A}$  in $\boldsymbol{\CC P^1 \times \CC P^1}$}

Consider the  anti-diagonal
\[
A:=\left\{ ([z:w],[\bar{z}:\bar{w}]) \, | \, [z:w] \in \CC P^1 \right\}
\]
which is a Lagrangian submanifold of $\CC P^1 \times \CC P^1$, where both factors of $\CP^1$ have the same standard symplectic form. Let $\mu : \CC P^1 \times \CC P^1 \to \RR^2$ be the moment map whose image is a square $P=[0, l]^2$ and $L_{(a,b)}$ be the moment f\/iber over $(a,b) \in P$.  Then, $L_{(a,b)}$ is a torus isomorphic to $L_1 \times L_2$ where $L_1$ has radius $a$ and $L_2$ has $b$. Note that if $a \neq b$,  $L_{(a,b)}$ does not intersect~$A$.
If $a=b$, they intersect along a circle.
\begin{figure}[t]
\centering
\includegraphics[height=4in]{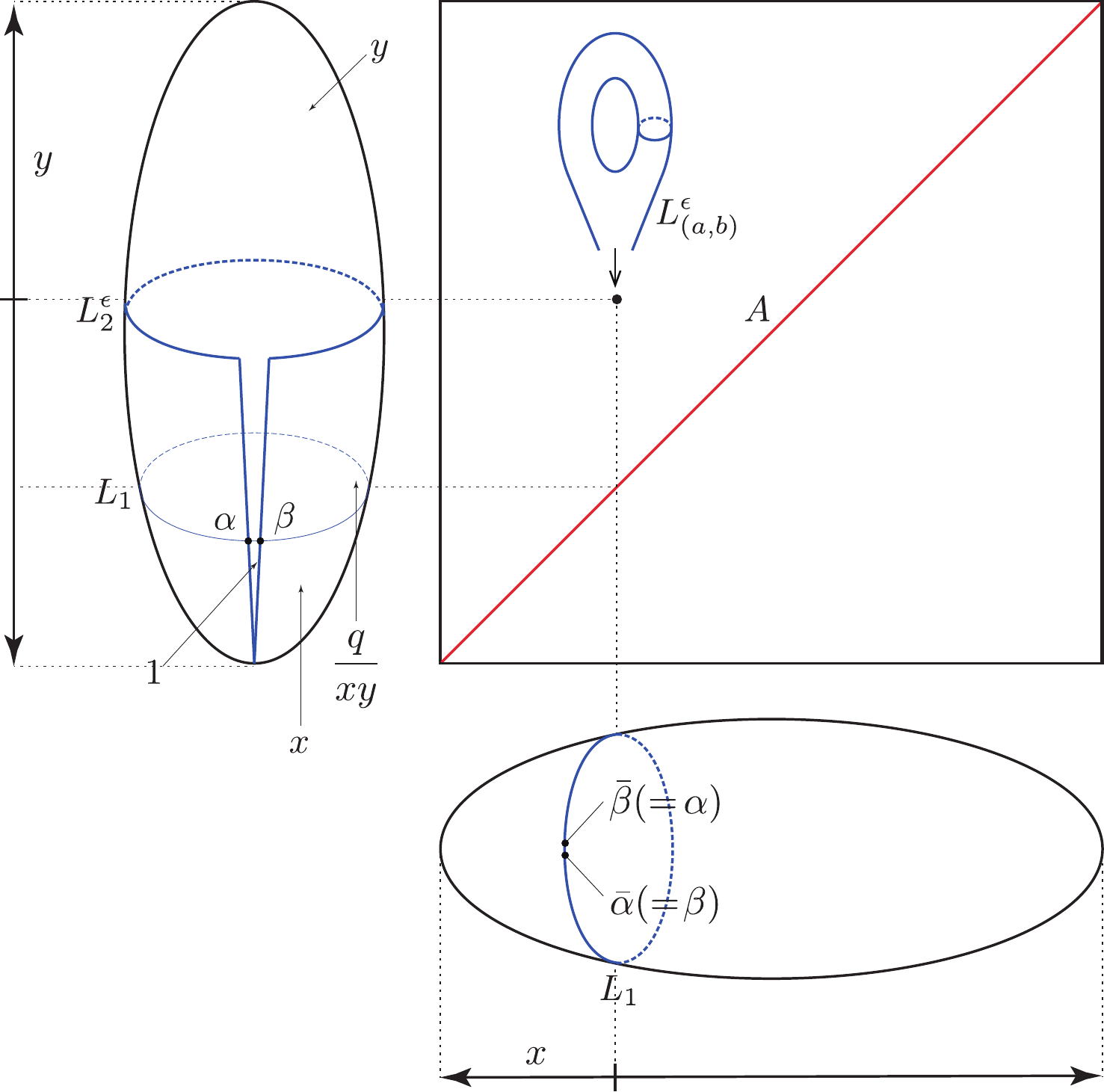}
\caption{Intersection of antidiagonal and the deformation $L^\epsilon_{(a,b)}$ of $L_{(a,b)}$.}
\label{fig:antidiag}
\end{figure}

In the example of $\CP^1$, the central f\/iber was deformed,
and its Floer chain complex with a~generic torus f\/iber was considered, whereas for our case of anti-diagonal,
we deform a generic torus f\/iber while keeping the anti-diagonal $A$ f\/ixed.
Namely, for each torus f\/iber, we deform the second component $L_2$ using the same methods as we did for $\CP^1$ in the previous section and get $L_2^{\epsilon}$ as in Fig.~\ref{fig:antidiag}.
In fact, consider the ``real'' circle in $\CP^1$ corresponding to the f\/ixed points of complex conjugation of $\CP^1$, and we
 may choose the deformation $L_2^\epsilon$ so that it  is symmetric with respect to this real circle.
 (In the Fig.~\ref{fig:antidiag}, the real circle is the vertical circle which bisects the spike of $L_2^{\epsilon}$.)

Then, $L_{(a,b)}^\epsilon :=  L_1 \times L_2^\epsilon$ will meet the anti-diagonal $A$ in at most two points. If $0 \leq a \leq b \leq l$,  they intersect precisely at two points and we can explicitly f\/ind out these two points. Let~$\alpha$ and~$\beta$ be two intersection points of~$L_1$ and $L_2^\epsilon$ as in the picture below. Then, it is easy to check that
\[
A \cap L_{(a,b)}^\epsilon = \{ (\alpha, \bar{\alpha}), (\beta, \bar{\beta}) \}.
\]
Note that $L_2^\epsilon$ will be preserved under the conjugation action on $\CC P^1$ and $(\alpha,\beta)$ and $(\beta,\alpha)$ are two intersection points of $L_1$ and $L_2^\epsilon$ since $\bar{\alpha} = \beta$ in this case

Now, we have to f\/ind holomorphic strips from $ (\alpha, \beta)$ to $ (\beta, \alpha) $ and vice versa. The following proposition classif\/ies all those strips in terms of holomorphic strips in the $\CC P^1$ (which one might think of as the f\/irst or the second factor of $\CC P^1 \times \CC P^1$).

We introduce the following notation. We say that a holomorphic strip $u:\R \times [0,1] \to M$ has a Lagrangian boundary condition $(L_a, L_b)$ if
the image of $(\R,0)$ maps to $L_a$ and that of $(\R,1)$ maps to $L_b$.

\begin{prop}
There is a one to one correspondence between holomorphic strips with bounda\-ry~$\big( A ,L^{\epsilon}_{(a,b)} \big)$ in $\CC P^1 \times \CC P^1$ and holomorphic strips with bounda\-ry~$(L_2^{\epsilon}, L_1)$ in $\CC P^1$. Moreover, corresponding strips have the same symplectic area.

The same holds for pairs $(L^\epsilon_{(a,b)}, A)$ and $(L_1, L_2^\epsilon)$.
\end{prop}

\begin{figure}[t]
\centering
\includegraphics[height=1.2in]{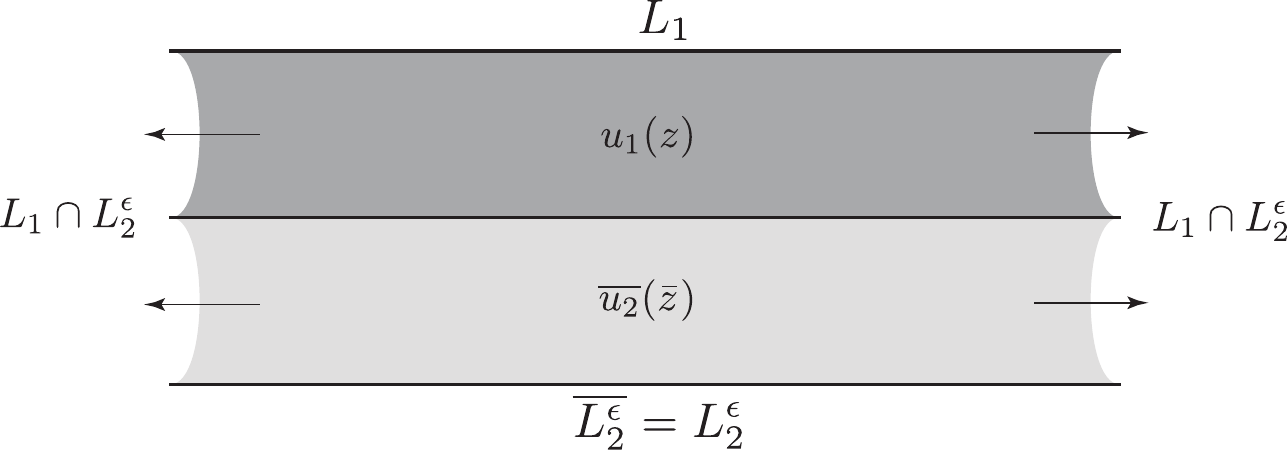}
\caption{Gluing of strips.}\label{fig5}
\end{figure}

\begin{proof}
Let $u = (u_1, u_2) : \RR \times [0,1] \to \CC P^1 \times \CC P^1$ be a holomorphic strip with boundary conditions
\[
u(\cdot,0 ) \in A , \qquad u(\cdot, 1) \in L_1 \times L_2^\epsilon,
\]
and with asymptotic conditions
\[
u(\infty,\cdot) = (\alpha, \beta), \qquad u(-\infty,\cdot) = (\beta, \alpha).
\]
From the boundary conditions of $u$, we can conclude that $u_1$ and $\overline{u_2}$ agrees on one of boundary components, i.e.\ $u_1 (s, 0) = \overline{u_2 (s,0)} $. Thus, if we def\/ine $u' : \RR \times [-1,1] \to \CC P^1 $ by (see Fig.~\ref{fig5})
\begin{gather*}
u' (z) =
\begin{cases}
  u_1 (z) = u_1 (s,t), & t \in [0,1],\\
 \overline{u_2} (\bar{z}) =  \overline{ u_2 ( s, -t)}, & t \in [-1,0],
\end{cases}
\end{gather*}
where we use the complex coordinate as $z = s + it$, then $u'$ asymptotes to $\alpha(=\bar{\beta})$ and $\beta(=\bar{\alpha})$ at $\infty$ and $-\infty$ respectively as we take complex conjugate of~$u_2$. Note that by the construction, $L_2^\epsilon$ is preserved by complex conjugation and hence, one of boundary components of the strip is still mapped to $L_2^\epsilon$ by~$\overline{u_2}$.

Finally, since $\CC P^1 \times \CC P^1$ has the product symplectic structure induced from one on each factor, $u'$ and $u$ should have the same symplectic area.
\end{proof}

Although the moduli spaces can be identif\/ied, the orientation of the moduli spaces work slightly dif\/ferently
(see~\cite{FOOO} for details on the def\/inition of canonical orientations).
For example, for a holomorphic strip with boundary $( L_2^\epsilon, L_1)$, changing the
orientation of~$L_2^\epsilon$ and $L_1$ to the opposite orientation reverses the canonical orientation of the holomorphic strip. But for a~holomorphic strip with boundary
$\big( A ,L^{\epsilon}_{(a,b)} \big)$, if we change the
orientation of $L_2^\epsilon$ and~$L_1$ at the same time, the orientation
of the product~$L^{\epsilon}_{(a,b)}$ remains the same, and so does the canonical orientation.
Hence, even though the holomorphic strips for the calculation of the Floer cohomology $HF(L_2^\epsilon, L_1)$ in $\CP^1$ cancels in pairs (to produce a non-vanishing Floer cohomology of the equator),
but the corresponding pairs of holomorphic strip with boundary $\big( A ,L^{\epsilon}_{(a,b)}\big)$
do not cancel because they have the same sign from this consideration. This is why all the terms in the matrix factorization below~\eqref{MFcp1cp1} has positive signs.

By the above proposition, to f\/ind holomorphic strips for the anti-diagonal, and a deformed generic torus f\/iber, it suf\/f\/ices to f\/ind holomorphic strips bounding~$L_1$ and $L_2^\epsilon$ which converge to~$\alpha$ and~$\beta$ at $\pm \infty$. These are the same holomorphic strips, discussed in the previous section for~$\CP^1$.  Namely, the shape of the strip remain the  same except that now we have deformed $L_2$ whose position is at $b$, not the central f\/iber of~$\CC P^1.$

Before we proceed, we recall the disc potential for~$\CP^1 \times \CP^1$,
whose terms correspond to holomorphic discs of index two, intersecting each toric divisor and having boundary lying in f\/ibers of the moment map
(see~\cite{FOOO2}):
\begin{gather*}
e^\alpha T^{u_1} + e^{-\alpha} T^{l-u_1} + e^\beta T^{u_2} + e^{-\beta} T^{l-u_2},
\end{gather*}
where $(\alpha,\beta)$ represents a induced holonomy on the boundary of holomorphic discs,
which may be identif\/ied with an element of $H^1$ of the torus. We denote $x = e^\alpha T^{u_1}$ and $y = e^{-\beta} T^{l-u_2}$ to obtain the potential (where $q = T^l$):
\[
W = x + \frac{q}{x} + y + \frac{q}{y}.
\]
\begin{Remark}
We use the coordinate $y = e^{-\beta} T^{l-u_2}$ instead of $y = e^{\beta} T^{u_2}$ so that the upper-hemisphere of the second factor $\CC P^1$ bounded by $L_2$ has the area $y$ (Fig.~\ref{fig:antidiag}). This is to get a~symmetric form of factorization of~$W$.
\end{Remark}

In Fig.~\ref{fig:antidiag}, strips from $\alpha$ to $\beta$ are strips of area~1 and $\frac{q}{xy}$, and those from $\beta$ to $\alpha$ are strips of area~$x$ and~$y$. Therefore, the resulting factorization of~$W$ is
\begin{gather}\label{MFcp1cp1}
(x+y)\left(1 + \dfrac{q}{xy} \right) =   x + \frac{q}{x} + y + \frac{q}{y}.
\end{gather}
These four strips contribute to $m_1$ with the same sign as we discussed above.
This proves Proposition~\ref{prop1}.

\begin{Remark}
We can also compute the matrix factorization corresponding to the central moment f\/iber. It turns out to be an exterior tensor product of matrix factorization of the central f\/iber of each factor $\CC P^1$ which is given in~\cite{CL}.
One can check that the following matrix factors $(W-\lambda) I_4$, where $\lambda= 4 \sqrt{q}$:{\samepage
\begin{gather*}
 \begin{pmatrix}0 & 0 & z-\sqrt{q} & -1+\frac{\sqrt{q}}{w} \vspace{1mm}\\
0 & 0 & w-\sqrt{q} & 1-\frac{\sqrt{q}}{z} \\
1-\frac{\sqrt{q}}{z} & 1-\frac{\sqrt{q}}{w} & 0 & 0 \\
-w+\sqrt{q} & z-\sqrt{q} & 0 & 0\end{pmatrix}.
\end{gather*}
For the tensor product of matrix factorization, see~\cite{ADD}.}
\end{Remark}

\section[Lagrangian Floer cohomology in $\CP^1 \times \CP^1$]{Lagrangian Floer cohomology in $\boldsymbol{\CP^1 \times \CP^1}$}\label{section5}

In this section, we verify the conjecture that in the derived Fukaya category,
the following two objects are the same:
\[
A \oplus A[1], \qquad T_{(1,-1)} \oplus T_{(-1,1)}.
\]
One is the direct sum of anti-diagonal $A$ and its
shift $A[1]$. The other is the direct sum of Lagrangian torus f\/iber at the center of
the moment map image with holonomy $(1,-1)$, denoted as $T_{1,-1}$ or $(-1,1)$, denoted
as $T_{-1,1}$. We denote by $T_0$ the central f\/iber of $\CP^1 \times \CP^1$
(without considering holonomies). We refer readers to Seidel's book \cite{S} on the def\/inition
of derived Fukaya category. We just recall that in our case, we
work with $\Z/2$-grading and by def\/inition we have $HF^*(L[1],L') = HF^{*+1}(L,L')$.

\subsection{Floer cohomology}

First we compute  $HF(T_{1,-1},A)$ and $HF(T_{-1,1},A)$.
Note that $T_0 \cap A$ is a clean intersection, which  is a circle $S^1$. Instead of working with the Bott--Morse version of the Floer cohomology, we move $T_0$ by Hamiltonian isotopy so that it intersects
$A$ transversely at two points. The Hamiltonian isotopy we choose are rotations in each factor of $\CP^1$
so that the equator of the circle is moved to the great circle passing through North and South pole.
More precisely, if we identify $\CP^1$ as $\C \cup \{ \infty \}$ and the equator with
the unit circle in~$\C$, then after isotopy, we obtain a Lagrangian submanifold~$L_0$ obtained
as a product of  real line in the f\/irst component and  imaginary line in the second component.

Locally on $\CC \times \CC$ we use $(a,b)$ and $(x,y)$ as coordinates of the f\/irst and the second factor, respectively. Let $L_0$ be the torus in $\CC P^1 \times \CC P^1$ given by the following equations
\begin{gather*}
L_0 =
\begin{cases}
a=0, \\
y=0.
\end{cases}
\end{gather*}
We denote $L^R$ (resp.~$L^I$) the great circle in $\CP^1$ corresponding to the real axis (resp. imaginary axis). We have $L_0 = L^R \times L^I $.

The anti-diagonal $A$ (for which we will write $L_1$ from now on) can be expressed as
\begin{gather*}
L_1 =
\begin{cases}
a-x=0, \\
b+y=0.
\end{cases}
\end{gather*}

Let us calculate the Floer cohomology of the pair $(L_0,L_1)$. They intersect at two points, $(0,0)$ and $(\infty, \infty)$.
We denote $p=(0,0)$ and $q=(\infty,\infty)$. As explained in the previous section,
given the holomorphic strips with boundary on $(L_0, L_1)$, we can glue the f\/irst and
the conjugate of second component of the strip to obtain a holomorphic strip with
boundary on $(L^I, L^R)$ (lower boundary on~$L^I$, and upper boundary on~$L^R$).

\begin{figure}[t]
\centering
\includegraphics[height=1.9in]{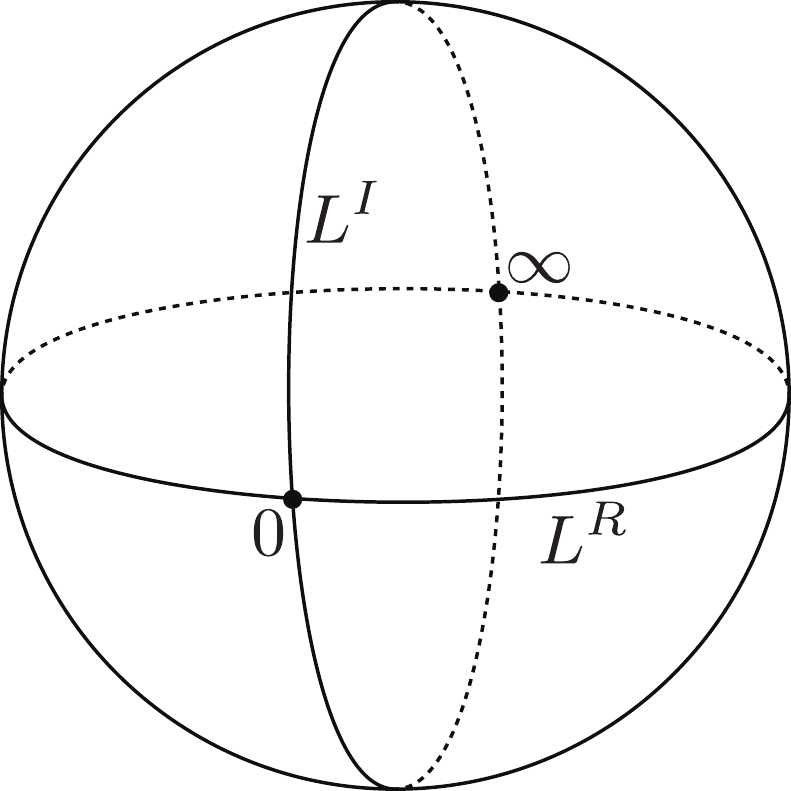}
\caption{The f\/irst factor of $\CC P^1 \times \CC P^1$.}
\label{rir}
\end{figure}

 There are four such strips as seen in Fig.~\ref{rir} (two strips from $p$ to $q$, the other two
 from~$q$ to~$p$) and these four strips have the same symplectic area.

As explained in the previous section, each of these strips are counted with the same sign in $\CP^1 \times \CP^1$
(dif\/ferent from the case of~$\CP^1$). Hence two strips from $p$ to $q$ do not cancel out but
adds up. In fact $m_1^2 \neq 0$ also, since
\[
m_1^2 =W_{T_0} - W_A = W_{T_0}
\]
and the potential $W_{T_0}$ for
the  central f\/iber $T_0$ with a trivial holonomy $(1,1)$ is non-trivial, which is
a sum of 4 terms corresponding to 4 holomorphic discs with  boundary either on $L^I$ or $L^R$
in~$\CP^1$.

But for $T_{1,-1}$ or $T_{-1,1}$ (which induces the f\/lat line bundles of the same holonomy on~$L_0$),
two strips from $p$ to $q$ cancel out due to holonomy contribution, and note that
the corresponding potential $W_{T_{1,-1}} = W_{T_{-1,1}} =0$. Thus from this cancellation of holomorphic strips we have
$m_1(p) =m_1(q)=0$.

Hence the Floer cohomology $HF(T_{1,-1}, A)$ (or $HF(T_{-1,1},A)$) is generated by $p$, $q$ and hence isomorphic  to the homology of $S^1$ with
Novikov ring coef\/f\/icient.
The similar argument works for $HF(A, T_{1,-1})$(or $HF(A, T_{-1,1})$), which is again generated by $p$, $q$.

The Floer cohomology $HF(T_{1,-1},T_{1,-1})$ is a Bott--Morse version of the Floer cohomology (see~\cite{FOOO}) and can be computed as in~\cite{C1} or~\cite{CO}, and is isomorphic to
the singular cohomology of the torus $H^*(T_0,\Lambda)$.
The Floer cohomology $HF(A,A)$ is also isomorphic to the singular cohomology $H^*(A,\Lambda)$,
as it is monotone and minimal Maslov index is 4 (see~\cite{Oh}).

\subsection{Products}
From now on, we don't distinguish $L_0$ and $T_0$ since they are clearly isomorphic in the Fukaya category. We assume that $L_0$ is the central torus f\/iber in $\CP^1\times \CP^1$ which is equipped with a~f\/lat complex line bundle of holonomy $(1,-1)$
(or $(-1,1)$), but we will omit it from the notation for simplicity. And by $L_1$, we denote the
anti-diagonal Lagrangian submanifold (for which we used the notation $A$ before).

\begin{Lemma}\label{lem1}
Consider the product
\[
m_2: \ HF(L_1, L_0) \times HF(L_0,L_1) \to HF(L_1,L_1)
\]
we have that $m_2(p,p) = [p] \pm[L_1] T^{l/2}$.
 Here $T^{l/2}$ is an area of
the upper $($or lower$)$ hemisphere of each factor $\CP^1$.
\end{Lemma}

\begin{Lemma}\label{lem15}
\[
m_2: \ HF(L_1, L_0) \times HF(L_0,L_1) \to HF(L_1,L_1)
\]
we have that $m_2(q,q) = [q]  \mp[L_1] T^{l/2}$.
\end{Lemma}

\begin{Lemma}\label{lem2}
\[
m_2: \ HF(L_0, L_1) \times HF(L_1,L_0) \to HF(L_0,L_0)
\]
we have that $m_2(p,p) = [p]\pm [L_0] T^{l/2}$.
\end{Lemma}

\begin{Lemma}\label{lem25}
\[
m_2: \ HF(L_0, L_1) \times HF(L_1,L_0) \to HF(L_0,L_0)
\]
we have that $m_2(q,q) = [q]\mp [L_0] T^{l/2}$.
\end{Lemma}
\begin{Lemma}\label{lem3}
For the product
\[
m_2: \ HF(L_1, L_0) \times HF(L_0,L_1) \to HF(L_1,L_1),
\]
 we have $m_2(p,q) = m_2(q,p) = 0$.
\end{Lemma}
\begin{Remark}
It turns out that the products $m_2(p,q)$, $m_2(q,p)$
for
\[
m_2: \ HF(L_0, L_1) \times HF(L_1,L_0) \to HF(L_0,L_0)
\]
do not vanish. But this won't be needed in our arguments of equivalence later
\end{Remark}

\begin{proof}[Proof of Lemma~\ref{lem1}.] The proof breaks into two parts, $(i)$~one for the actual counting of strips and $(ii)$~the other for the Fredholm regularity of these strips.

$(i)$ The holomorphic triangle contributing to $m_2$ in this case can be considered as
a holomorphic strip $u : \R \times [0,1] \to \CP^1\times \CP^1$
with
\[
u(\cdot, 0) \subset L_0, \qquad  u(\cdot,1) \subset L_1
\]
and a marked point $z_0 = [t_0,1]$ on the upper boundary of the strip, which is used
as an evaluation to $L_1$.
Hence the f\/irst and second factor of holomorphic triangle can be again glued as in the previous
section to give a holomorphic strip in $\CP^1$ with boundary on
$(L^R, L^I)$ in  $\CP^1$,  but both ends of the holomorphic strips converge to~$0$
(the f\/irst component of~$p$). For convenience, we also call~$0$ as~$p$.

Note that both $L^R$ and $L^I$ are preserved by the complex conjugation so that we can freely use this kind of process. Note also that after gluing, the marked point for evaluation lies in the interior of the strip.

From \cite{KO} such holomorphic strips can be decomposed into simple ones, and
in this case, homotopy class of any holomorphic strip is given by the union of strips (in fact an even number of unions to come back to~$p$).
Since we are interested in the case that  the dimension of the evaluation image is either 0 or two,
the number of strips must be less than or equal to two. Since it starts and ends at~$p$, the number is either 0 or 2.

First we consider the case of 0, or equivalently a constant triangle.
In this case, we can use the following theorem of the f\/irst author in preparation

\begin{Theorem}[\cite{C2}]
Let $L_a$, $L_b$, $L_c$ be Lagrangian submanifolds in a $2n$-dimensional symplectic manifold $M$, such that all possible intersections among them are clean. If $L_a \cap L_b \cap L_c =\{p\}$,
$p$~contributes to energy zero part of the product  $([L_a \cap L_b])  \times ([L_b \cap L_c])$ in $HF(L_a,L_b) \times HF(L_b, L_c)$ non-trivially as $p=L_a\cap L_b \cap L_c$,  if and only if
\[
\dim_\R(L_a \cap L_b) + \dim_\R(L_b \cap L_c) + \dim_\R(L_c \cap L_a)+ \angle L_aL_b + \angle L_bL_c + \angle L_cL_a=2n.
\]
\end{Theorem}
The notion of an angle is def\/ined in~\cite{Al}. In our case, $L_a=L_c= L_0$ and $L_b=L_1$ and
hence $\angle L_cL_a=0$ and  also it is not hard to see from the def\/inition of an angle that
if $L_a$ and $L_b$ intersect transversely,
\[
\angle L_aL_b + \angle L_bL_a = n.
\]
Thus, $n + \angle L_aL_b + \angle L_bL_c + \angle L_cL_a$  is  equal to $2n$.
Hence, we have $p $ as an energy-zero component of $m_2(p,p)$.

Now, we consider the case that the holomorphic strip covers half of~$\CP^1$. It is in fact easy to f\/ind such holomorphic strips, covering half of $\CP^1$,
starting from $p$ ending at $p$.  The image of the strip is a disc with boundary either on real or
imaginary circle in $\CP^1$, and  one of the lower or upper boundary covers the circle once, and
the other boundary covers part of the segment and comes back to~$p$. (This strip of index two usually appears to explain the bubbling of\/f in $\C$ with Lagrangian submanifolds $\R$ and unit circle~$S^1$.)

\begin{figure}[t]
\centering
\includegraphics[height=1.5in]{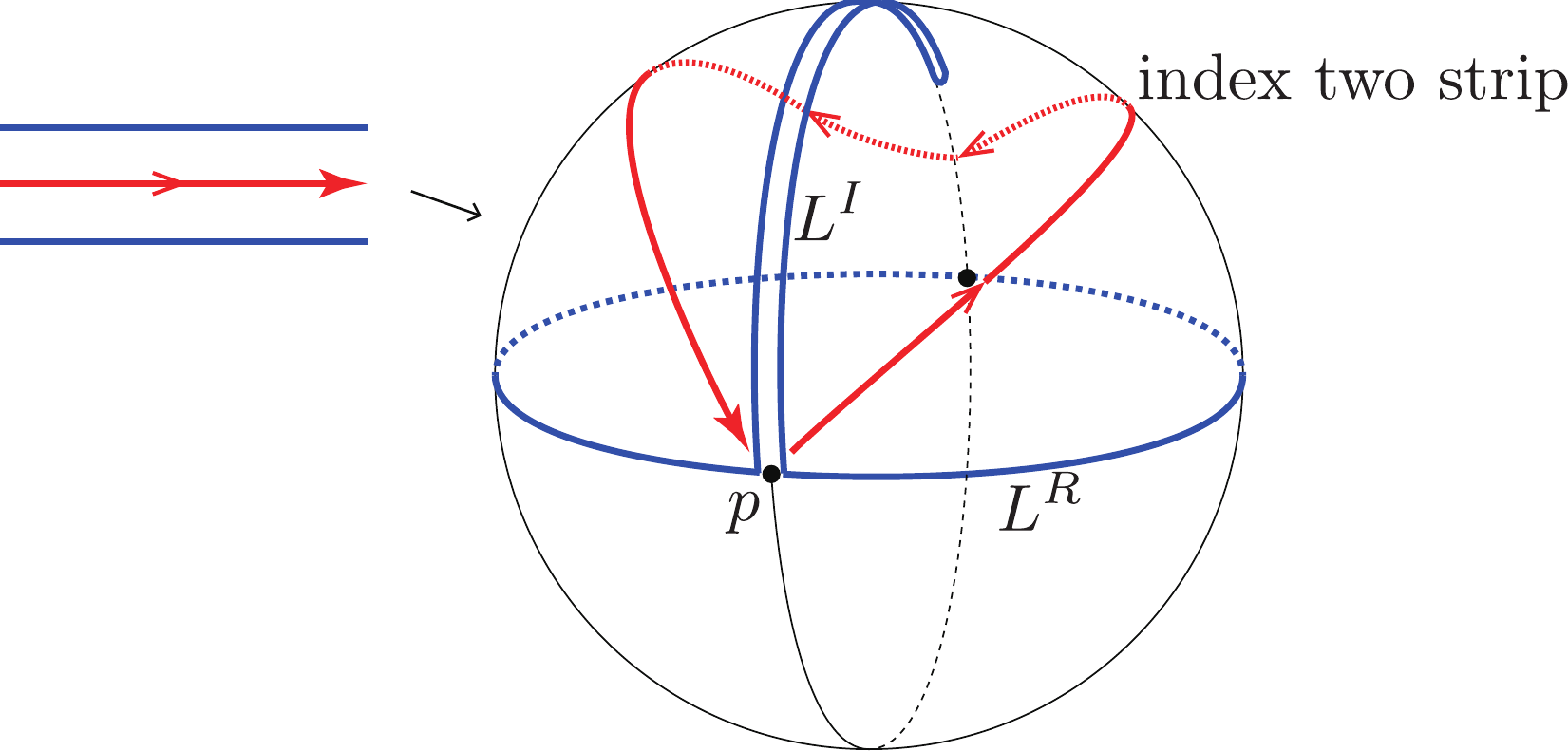}
\caption{Index two strips bounding $L^I$ and $L^R$.}\label{fig7}
\end{figure}

These holomorphic strips of Maslov--Viterbo index two, lies inside a holomorphic disc in $\CP^1$
(of Maslov index two) with boundary on $L^I$ or $L^R$, and there are 4 such discs (see Fig.~\ref{fig7}).
Thus, there are 4 homotopy classes of  Maslov--Viterbo index two holomorphic strips from~$p$ to~$p$ and  we denote them as $\beta_1,\dots ,\beta_4$.

Consider  $\CM_1(L_0,L_1,\beta_i, p,p)$ the moduli space of holomorphic strips of  (Maslov--Viterbo) index two as described above starting and ending at~$p$, with one marked point in the upper boundary of
the strip for $i=1,\dots, 4$. The boundary $\partial \CM_1(L_0,L_1,\beta_i,p,p)$ is
well understood, and exactly has two possible components, one is from the broken strip
of from~$p$ to~$q$ and to~$p$, and the other is the bubbling of\/f of a Maslov index two disc attached
to a constant strip at $p$. In the former case, the marked point is located in either component of the broken holomorphic strips, and in the latter case, one of the coordinate of the marked point is free to move along the bubbled disc.

We can compare the orientations of the bubbled discs for each $\beta_i$'s and they correspond to
the potential $W$ of $L_0$, and with the holonomy $(1,-1)$ or $(-1,1)$, all these terms cancel out.
Similarly, the evaluation images of the f\/irst type of boundary from the broken strips also are mapped  exactly twice, since given an index one strip, there are two adjacent strips to it. And  as the
signs cancelled in~$W$, the signs of the images for the f\/irst type of boundary should be opposite too.
Thus, this shows that actually the boundaries of $\CM_1(L_0,L_1,\beta_i,p,p)$ for $i=1,\dots,4$
matches with opposite signs and the union gives a cycle in $L_1 \cong \CP^1$.

\begin{figure}[t]
\centering
\includegraphics[height=1.5in]{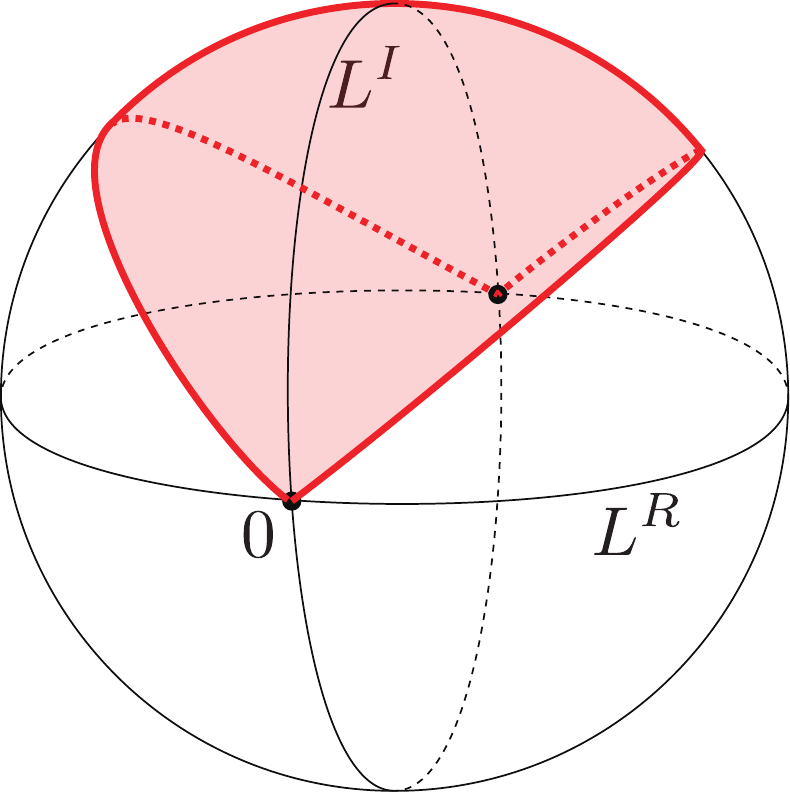}
\caption{The shape of a lune.}
\label{lune}
\end{figure}

Hence this shows that $m_2(p,p)$ is a constant multiple of $[L_1]$. And it is enough to
f\/ind the constant.
Given such an index two strip, we consider the glued strip in $\CP^1$, and  we evaluate at the marked point which is in the middle line of the glued strip. By varying the strip, it is not hard to
see that  the image of evaluation map covers ``half'' of the disc, or a spherical lune (Fig.~\ref{lune}), connecting~$p$ and~$q$ once. But there are 4 discs and these 4 lunes together cover the whole $\CP^1$. This shows that
the constant is one, and we have $m_2(p,p) = \pm [L_1] T^{l/2}$.

$(ii)$ One can show that these strips are Fredholm regular from the following explicit formulation.
First we identify the holomorphic strip with the upper half-disc $D_{+} =\{z \in \C \mid |z| \leq 1$, ${\rm Im}\,(z) \geq 0\}$
with punctures at $-1,+1 \in D^2$, which are identif\/ied with  $-\infty, \infty$ of the strip.
Then, consider a  holomorphic map $u : D_+  \to \C$ with semi-circle of $\partial D_+$ mapping to the unit circle
of~$\C$, and real line segment of $\partial D_+$ mapping to real line of~$\C$. All such maps of
degree two (whose images covers~$D^2$ once) are given by
\begin{gather}\label{eq1}
z \in D_+ \mapsto  \frac{(z-a)(z-b)}{(1-az)(1-bz)}
\end{gather}
for a real number $a, b \in (-1,1)$,
or
\begin{gather}\label{eq2}
z \in D_+ \mapsto  \frac{(z-\alpha)(z-\overline{\alpha})}{(1-\alpha z)(1- \overline{\alpha}z)}
\end{gather}
for some $\alpha \in D^2$.

(To see this one starts with the generic form of a product of two Blaschke factors, and def\/ine
an involution $u(z) \to \overline{u(\overline{z})}$ and f\/ind its f\/ixed elements.)
 Since the holomorphic discs in $\C$ with boundary on $S^1$  are always Fredholm regular,
the f\/ixed elements by involution are again Fredholm regular.
\end{proof}

\begin{proof}[Proof of Lemma~\ref{lem15}.]
All the arguments are the same as the proof of Lemma~\ref{lem1} except on the sign in front of~$[L_1]$. Hence, we only need to compare orientations. Note that the moduli space of holomorphic strips from $p$ to $p$ of index two with boundary on $(L^R, L^I)$ in $\CP^1$, gives rise to the moduli space of holomorphic strips from $q$ to $q$ by rotating $180$ with boundary on $(L^R, L^I)$ in~$\CP^1$. If $L^R$ lies at the center of the disc which contains the strip, then
this  this process reverses the orientation of $L^R$, but not the orientation of $L^I$.
(If $L^I$ lies at the center, orientation of $L^I$ is reversed, but that of $L^R$ is f\/ixed.)

As the rest of the ingredients for the orientation of the moduli space and evaluation map to the anti-diagonal remain the same, the resulting evaluation image has the opposite sign.
\end{proof}

\begin{proof}[Proof of Lemma~\ref{lem2}.]
The proof is somewhat similar to that of Lemma~\ref{lem1}.

A holomorphic triangle contributing to $m_2$ in this case can be considered as
a holomorphic strip $u : \R \times [0,1] \to \CP^1\times \CP^1$
with
\[
u(\cdot, 0) \subset L_1, \qquad  u(\cdot,1) \subset L_0
\]
and a marked point $z_0 = [t_0,1]$ on the upper boundary of the strip, which is used
as an evaluation to $L_0$.

Hence the f\/irst and second factor of holomorphic triangle can be again glued as in the previous
section to give a holomorphic strip in $\CP^1$ with boundary on $(L^I, L^R)$ of $\CP^1$, but
with both ends converging to~$p$.

Again, the same argument as in the previous lemma shows that constant strip do contribute to~$m_2$ in this case, and also the Maslov--Viterbo index two strips are to be considered. The relevant
moduli space of holomorphic strips have 4 connected components, and their boundaries cancel out.
Hence the evaluation image def\/ines a 2-dimensional cycle in $L_0$, or a constant multiple of unit $[L_0]$.

Thus it is enough to f\/ind the constant. For this we use the explicit form of the holomorphic strip
\eqref{eq1}, \eqref{eq2}. We may f\/ind a holomorphic strip
of index two, sending $z_0$ to $(t_1,t_2) \in L_0$.
After gluing of f\/irst and second component of the strip, we may f\/ind a holomorphic map
from~$D_+$ sending $0$ to $t_1 \in \R \subset \C$ and~$i$ to $\overline{t_2} \in S^1 \subset \C$
(up to automorphism of a strip, we may assume that $(0,i)$ corresponds to $z_0$).

By inserting these numbers to \eqref{eq1}, we obtain
\begin{gather*}
a + b  =  \frac{(t_1-1)(\overline{t_2} +1)}{i (1-\overline{t_2})}, \qquad
ab   =  t_1.
\end{gather*}
Or from \eqref{eq2}, we obtain
\begin{gather*}
\alpha + \overline{\alpha}  =  \frac{(t_1-1)(\overline{t_2} +1)}{i (1-\overline{t_2})},\qquad
\alpha \overline{\alpha}   =  t_1.
\end{gather*}
Thus $a$, $b$ or $\alpha$, $\overline{\alpha}$ are (real or conjugate) pair of solutions of
the quadratic equation
\[
x^2  -\frac{(t_1-1)(\overline{t_2} +1)}{i (1-\overline{t_2})} x + t_1 =0.
\]
(one can check that the coef\/f\/icient of $x$ is real).
If we choose $t_1 < 1$ to be almost as big as~$1$,
and choose $t_2$ to be close to~$-1$, then the coef\/f\/icient of~$x$ is very close to~$0$ whereas~$t_1$ is almost equal to~$1$. Thus the quadratic equation has a unique conjugate pair of complex
solutions, both of which  lies in the unit disc (since $|\alpha|^2 <1$).
Thus, this shows that the constant~$c$ of~$m_2(p,p) = c[L_0]T^{l/2}$ equals~$\pm 1$.
Hence, this proves the lemma.
\end{proof}

The proof of Lemma~\ref{lem25} is exactly the same as that of Lemma~\ref{lem15}
and omitted.

\begin{proof}[Proof of Lemma~\ref{lem3}.]
We begin the proof of Lemma~\ref{lem3}.
The products
\[
HF(L_1, L_0) \times HF(L_0,L_1) \to HF(L_1,L_1),
\]
given by $m_2(p,q)$ or $m_2(q,p)$ are zero since
$HF(L_1,L_1) \cong H^*(\CP^1,\Lambda)$ has no degree one classes.
(This is because holomorphic strips connecting~$p$ and~$q$ have
odd  Maslov--Viterbo index, which is the dimension of the moduli space of holomorphic
strips).
\end{proof}

\subsection{Floer cohomology between torus with dif\/ferent holonomies}

First we consider the case of a cotangent bundle of a torus.
Let $L$ be a Lagrangian torus $T^n \subset T^*T^n$
Let $\CL_1$ and $\CL_2$ be two dif\/ferent f\/lat line bundles on~$L$.
We prove that the Floer cohomology of the pairs
$HF\big( (L,\CL_1), (L,\CL_2))$ vanishes if $\CL_1 \neq \CL_2$.

\begin{prop}
The Floer cohomology $HF((L,\CL_1),(\phi(L),\phi_*(\CL_2))$ vanishes
if $\CL_1 \neq \CL_2$.
\end{prop}
\begin{proof}
 Since $L$ is a torus, we identify $L$ as $\RR^n/ \ZZ^n$, and def\/ine a Morse
function $f:T^n \to \RR$ by
\begin{gather}\label{eq:morsett}
f(x_1,x_2,\dots,x_n) = \sum_{i=1}^n \cos (2\pi x_i).
\end{gather}
It is immediate to check that the critical points set is
\[
\{ (a_1,a_2, \dots,a_n) | \ a_i = 0  \ \textrm{or} \ 1/2 \  \textrm{for} \ i = 1,\dots,n \}.
\]
Denote the holonomy of $\CL^0$ (or $\CL^0$) along the $i$-th generator of $T^n$
by $h_i^0$ (or $h_i^1$).
Let $I = (1/6,5/6), J=(4/6,8/6)$.
Def\/ine
\[
\mathcal{S} = \{ L_1\times \cdots \times L_n \subset \RR^n/ \ZZ^n\,\,|
\  L_i = I \ \textrm{or} \  J \  \textrm{for} \ i = 1,\dots,n \}.
\]
$\mathcal{S}$ def\/ines an open cover of $T^n$. The line bundle $\CL^0$ (and $\CL^1$)
may be
described by local charts on $\mathcal{S}$ as follows.
We explain how to glue trivial lines bundles on the open sets of $\mathcal{S}$
\[
\phi: \ L_1\times \cdots \times L_n \times \CC \mapsto
L_1'\times \cdots \times L_n' \times \CC
\]
sends $(x_1,\dots,x_n,l) \to (x_1,\dots,x_n,l')$
where $l' = b_1 b_2 \cdots b_n l$ with
\[ b_i = \begin{cases}
 1 & \textrm{if} \ L_i = L_i', \\
 1 & \textrm{if} \ x_i \in (1/6,2/6), \\
 h_i^0 & \textrm{if} \ x_i \in (4/6,5/6) \ \textrm{and} \ L_i = I,\ \ L_i' = J, \\
 1/h_i^0 & \textrm{if} \ x_i \in (4/6,5/6) \ \textrm{and} \ L_i = J,\ \ L_i' = I.
\end{cases}
\]
It is easy to check that this def\/ines the f\/lat line bundle $\CL^0$.

Now we compute the boundary map in the Floer complex.
First, we f\/ix some sign convention about Morse complex.
Recall the following rules, for a submanifold $P \subset L$ and $x \in P$,
\[
 N_x P \oplus T_x P = T_x L.
 \]
Also
\[
 N_x P_1 \oplus N_x P_2 \oplus T_x(P_1 \cap P_2) = T_x L
 \]
determines the orientation of $P_1 \cap P_2$ at $x$.
Now, we denote $W^u(x)$, $W^s(x)$ to be the unstable and stable manifold
of $x$ for the given Morse function $f$ on $L$.
Then, we set
\begin{gather}\label{eq:morses1}
TW^{s}(x) \oplus TW^{u}(x) = T_x L.
\end{gather}

Finally, we set the orientation of the moduli space $\CM(x,y)$ of
the trajectory moduli space as
\[
 W^s(y) \cap W^u(x) = \CM(x,y).
 \]

Now, we consider the function $f$ given by (\ref{eq:morsett}).
Unstable manifolds of $f$ can be written as products of
intervals $[0,1/2)$ or $(1/2,0]$, and intervals are
canonically oriented. Hence we assign the
product orientations on the unstable manifolds.

\begin{Lemma}
Let $x=[a_1,a_2,\dots,a_n]$, $y=[b_1,\dots,b_n]$
where for a fixed $i$, $a_i=0$, $b_i = 1/2$ and $b_j= a_j$ for $j \neq i$.
Then, the trajectory space $\CM(x,y)$ has the canonical orientation
$(-1)^A \partial_i$ where $A$ is the number of $j <i$ with $a_j =0$.
Here $\partial_i$ is the $i$th standard basis vector of $\RR^n$.
\end{Lemma}
\begin{proof}
First, from the orientation convention, we can identify
$NW^{u} = TW^{s}$. Hence,
\begin{gather}\label{eq:morses2}
NW^{s}(y) \oplus NW^{u}(x) \oplus T\CM(x,y)
= NW^{s}(y) \oplus TW^{s}(x) \oplus T\CM(x,y) = TL.
\end{gather}
It is easy to check that
\[
(-1)^A \partial_i \oplus TW^{u}(y) = TW^{u}(x),
\]
where $A$ is the number of $j <i$ with $a_j =0$,
by comparing two unstable manifolds.
Hence, from~(\ref{eq:morses1}), we have
\[
 TW^{s}(x) \oplus (-1)^A \partial_i  = TW^{s}(y).
 \]
Hence, combining with (\ref{eq:morses2})
and denoting $T\CM(x,y) = (-1)^B\partial_i$, we have
\[
NW^{s}(y) \oplus TW^{s}(y) \cdot (-1)^A (-1)^B = TL.
\]
Hence, we have $TL \cdot (-1)^{A+B} = TL$, which proves
the lemma.
\end{proof}

The lemma implies that actual Morse boundary map
is given as follows by comparing the coherent orientation
with the f\/low orientation.
\[
\partial_{\rm Morse} x = (-1)^A (1 - 1) y = 0.
\]

Now, in the case of the Floer complex twisted by f\/lat bundles, we have
\begin{gather*}
\partial\big((a_1,\dots,a_n)\big) =  \sum_{\textrm{for each} \ a_i=0} (-1)^{A_i}
 =\left(1-\dfrac{h_i^0}{h_i^1} \right)(a_1,\dots,a_{i-1},a_i+1/2,a_{i+1},\dots,a_n).
\end{gather*}

If $\CL^0 = \CL^1$, we have $h_i^{0}/ h_i^{1} = 1$, hence  all boundary maps vanish and we obtain the singular cohomology of the torus $T^n$.
If $\CL^0 \neq \CL^1$,  we f\/irst
assume that $h_i^0 \neq h_i^1$ for all $i$, and show that
the complex has vanishing homology.

In fact, the above complex, with an assumption $h_i^0 \neq h_i^1$ for all $i$,
is chain isomorphic to the same complex with the following new dif\/ferential
\[
\widetilde{\partial}\big((a_1,\dots,a_n)\big) = \sum_{\textrm{for each} \ a_i=0} (-1)^{A_i} (1)(a_1,\dots,a_{i-1},a_i+1/2,a_{i+1},\dots,a_n).
\]
Here chain isomorphism can be def\/ined as
\[
\Psi([a_1,\dots,a_n]) = \left( \prod_{i \; \textrm{with} \ a_i =0}
(1-h_i^{0}/ h_i^{1}) \right) [a_1,\dots,a_n].
\]
It is easy to check that $\Psi \partial = \widetilde{\partial} \Psi$,
and there is an obvious inverse map.

The new complex with $\widetilde{\partial}$ may be considered
as the reduced homology complex of the standard simplex
$\Delta^{n-1}$, hence has a vanishing homology.
The face corresponding to $[a_1,\dots,a_n]$
contains $i$-th vertex if and only if $a_i =0$.

Now, consider the general case that $h_i^{0} \neq h_i^1$ if and only if $i \in \{i_1,\dots, i_k\}$
where $k \geq 1$.
The chain complex we obtain has non-trivial dif\/ferential only for the terms
containing $(1-h^0_i/h^1_i)$ for $i \in \{i_1,\dots, i_k\}$ and hence the chain complex decomposes
into several chain sub-complexes with only non-trivial dif\/ferentials within. And
by using the result in the f\/irst case, we obtain the proposition.
\end{proof}

So far, we have discussed the case in the cotangent bundle of the torus.
For our case,  it follows from the spectral sequence of~\cite{Oh2}.
\begin{Lemma}\label{lem4}
\[
HF(T_{1,-1},T_{-1,1}) \cong 0.
\]
\end{Lemma}
\begin{proof}
$T_0$ (and hence $T_{1,-1}$ and $T_{-1,1}$) is a monotone Lagrangian submanifold, and hence by~\cite{Oh2}, there is f\/iltration of the
Floer dif\/ferential
\[
m_1 = m_{1,0} + m_{1,N} + m_{1,2N} + \cdots,
\]
where $N$ is the minimal Maslov number of~$T$, which can be
easily modif\/ied to the case of f\/lat complex line bundles. By usual  spectral sequence argument,
we obtain the vanishing of the homology of~$m_1$ dif\/ferential, since the homology of
$m_{1,0}$ vanishes in our case.
\end{proof}

\subsection{Equivalence}

To show that $A \oplus A[1]$ is  equivalent  to $T_{1,-1} \oplus T_{-1,1}$.
We f\/ind
\begin{gather*}
\Phi_1 \in \Hom_{\rm DFuk}(A \oplus A[1], T_{1,-1} \oplus T_{-1,1}), \qquad
\Phi_2 \in \Hom_{\rm DFuk}(T_{1,-1} \oplus T_{-1,1},A \oplus A[1] ),
\end{gather*}
such that $\Phi_1 \circ \Phi_2 = {\rm Id}$, and $\Phi_2 \circ \Phi_1 = {\rm Id}$ in the derived Fukaya category of $\CP^1 \times \CP^1$.
We write
\begin{gather*}
\Phi_i =  \begin{pmatrix}\alpha_i & \beta_i \\\gamma_i & \delta_i\end{pmatrix}.
\end{gather*}
Namely,
\begin{gather*}
\alpha_1 \in \Hom(A, T_{1,-1}), \beta_1 \in \Hom(A, T_{-1,1}),\\
 \gamma_1 \in \Hom(A[1], T_{1,-1}), \delta_1 \in \Hom(A[1], T_{-1,1}),\\
 \alpha_2 \in \Hom(T_{1,-1}, A), \beta_2 \in \Hom( T_{-1,1},A[1]),\\
 \gamma_2 \in \Hom(T_{-1,1}, A), \delta_2 \in \Hom(T_{-1,1},A[1]).
\end{gather*}

We choose \begin{gather*}
\Phi_1 =  \begin{pmatrix} p & q \\ q &  p \end{pmatrix},\qquad
\Phi_2 = \frac{1}{2T^{l/2}} \begin{pmatrix}  p & -q \\ -q & p \end{pmatrix}.
\end{gather*}

\begin{Theorem}\label{thm:equi}
We have
\begin{gather*}
\Phi_1 \circ \Phi_2 = \pm {\rm Id} \in \Hom(A \oplus A[1], A \oplus A[1]),\\
 \Phi_2 \circ \Phi_1 = \pm {\rm Id} \in \Hom(T_{1,-1} \oplus T_{-1,1}, T_{1,-1} \oplus T_{-1,1}).
 \end{gather*}
 Therefore in the derived Fukaya category of $\CP^1 \times \CP^1$, $A \oplus A[1]$ is  equivalent  to $T_{1,-1} \oplus T_{-1,1}$.
\end{Theorem}

\begin{proof}
This follows from Lemmas \ref{lem1}, \ref{lem15}, \ref{lem2}, \ref{lem25} and \ref{lem3} and
the fact that $[p]=[q]$ in the Bott--Morse Floer cohomology of $A$ or $T$.
We note that for $\Phi_2 \circ \Phi_1$, the
products of the following type,
\[
HF(T_{1,-1}, A) \times HF(A, T_{-1,1}) \to HF(T_{1,-1}, T_{-1,1}),
\]
are automatically zero, due to Lemma~\ref{lem4}.

Note also that $[L_0]$ and $[L_1]$ play a role of units in $HF(L_0,L_0)$ and $HF(L_1,L_1)$, respecti\-vely.
\end{proof}

\section{Teardrop orbifold}\label{SecTO}

We show that the correspondence between the Floer complex equation and the matrix facto\-ri\-zation continues
to hold for a teardrop orbifold. Such correspondence can be divided into two levels.
The f\/irst level is regarding smooth discs. Namely, we consider the Floer complex equation, only involving
smooth holomorphic strips (and discs). Then, we obtain a smooth potential or the Hori--Vafa
Landau--Ginzburg potential and the correspondence holds on this level. Here, by smooth holomorphic strips or discs, we mean a holomorphic maps from a smooth domain Riemann surface with boundary, and
by def\/inition of holomorphicity,  they locally lift to uniformizing chart of the target orbifold point, and hence when their images contain an orbifold point, it meets the point with multiplicity (see~\cite{CP}).

For the second level
  we  consider bulk deformations by twisted sectors,
and hence obtain the corresponding bulk potential or bulk Landau--Ginzburg potential which
has additional terms corresponding to orbi-discs. We consider the
Floer complex equation, involving smooth and orbifold holomorphic strips (and discs), which
are maps from orbifold Riemann surfaces with boundary. Then the
correspondence  between the Floer complex equation and the matrix factorization continues
to hold for bulk deformed cases.

Let $X$ be the orbifold obtained from the following stacky fan:
\begin{figure}[h]
\centering
\includegraphics[height=1.3in]{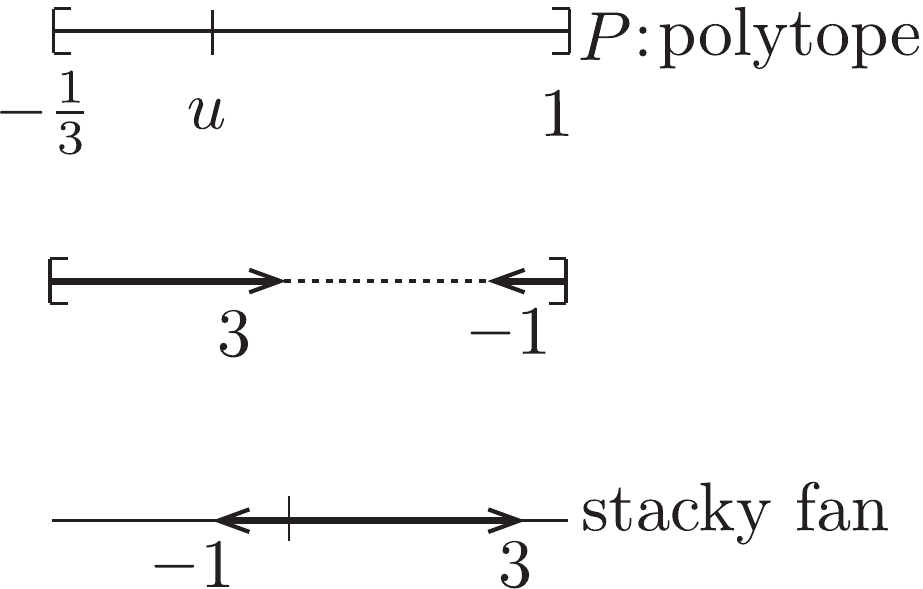}
\end{figure}
Then, $X$ is an orbifold with one singular point with $(\ZZ / 3\ZZ)$-singularity.

\subsection[The case of Hori-Vafa potential]{The case of Hori--Vafa potential}
The Hori--Vafa Landau--Ginzburg potential can be constructed (see~\cite{CP}) in this case by
consi\-de\-ring smooth holomorphic discs of Maslov index two:
\begin{gather}\label{teardroppoten}
W(z) = z ^3+ \frac{q^4}{z},
\end{gather}
where $z^3$ is due to the fact that smooth holomorphic discs around the orbifold point has to wrap
around it 3 times (see \cite{FOOO2} for a general procedure
of boundary deformation to construct such a potential from the moduli of holomorphic discs).

We brief\/ly review how to obtain the above expression of the potential~\eqref{teardroppoten} as above. Since index 2 discs correspond to the vectors in the stacky fan~\cite{CP}, we have the following description of index 2 holomorphic discs,
\[
e^{3x} T^{3 \left(u -  \left(-  \frac{1}{3} \right) \right)} + e^{-x} T^{ 1 -u} = e^{3x} T^{1+3u} + e^{-x} T^{1-u},
\]
where the power of $e$ represents the holonomy factors, and that of $T$ represents the area of discs (see~\eqref{eq:matfun}). In particular, $u$ is a position in the interior of the moment polytope. Note that we multiply $3$ to $(u+ 1/3)$ to obtain the area of the {\it smooth discs}.

One get the expression of $W$ as in \eqref{teardroppoten}, by substituting $z=e^x T^{1+3u}$ and $q=T^{1/3}$. Then, the total area of the teardrop orbifold will give the term
\[
T^{1- (-1/3)} = T^{4/3} = q^4.
\]

Denote $L_u$ by the torus f\/iber over $u \in [-1/3,1]$ where we identify $P$ with the interval $[-1/3,1] \subset \RR$. Let $L$ be the balanced f\/iber $L_0$ (i.e. the moment f\/iber over $u=0$).

\begin{figure}[t]
\centering
\includegraphics[height=3.3in]{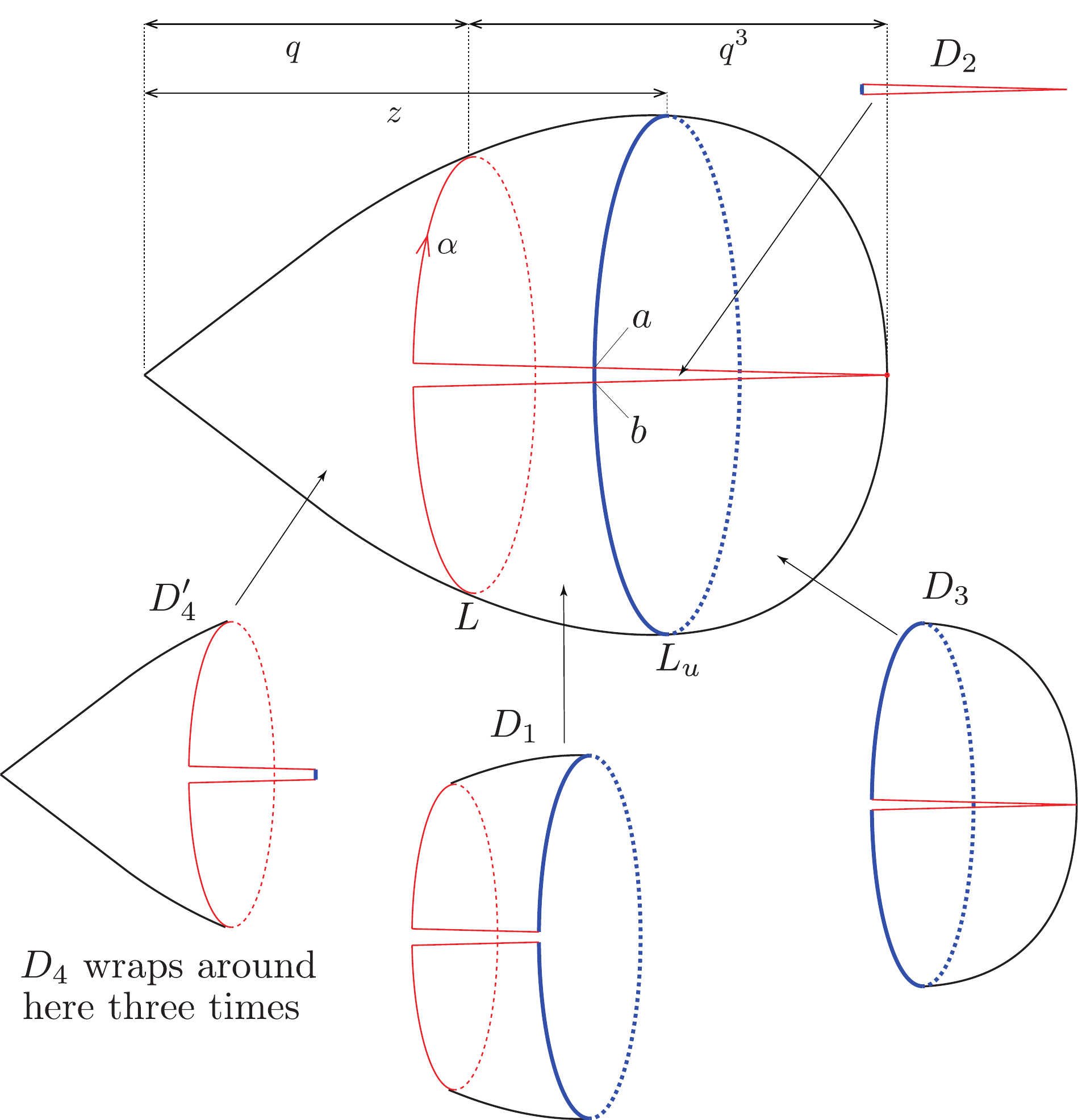}
\caption{Disk splitting in the teardrop orbifold.}
\label{dsteardrop}
\end{figure}

\begin{figure}[tp!]
\centering
\includegraphics[height=3.8in]{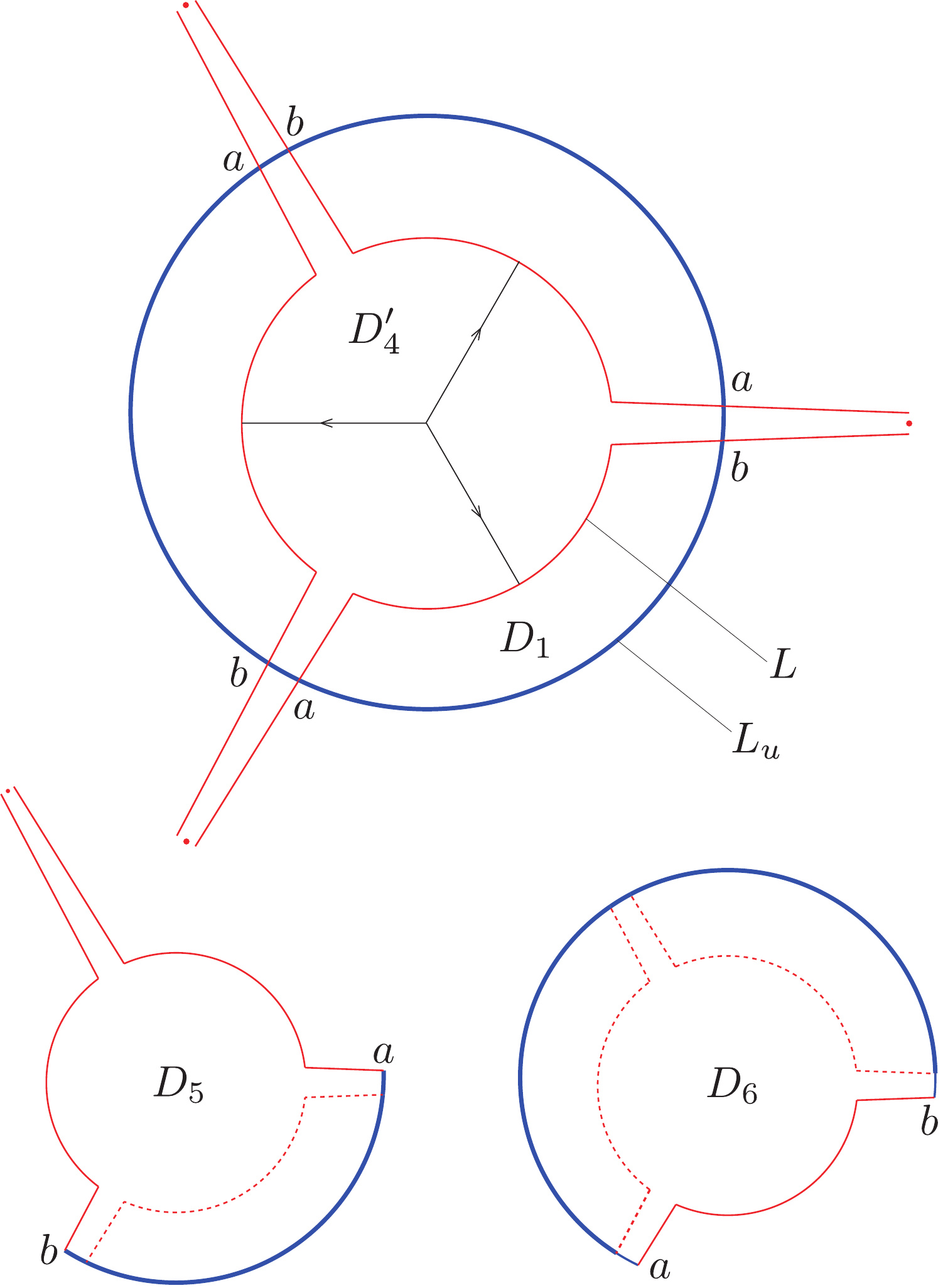}
\caption{Development f\/igure.}\label{develop}
\end{figure}

Let $\alpha$ be a holonomy around $L$ which is one of solutions of
\[
3 z^2 - \dfrac{1}{z^2} =0
\]
or equivalently, $3 \alpha^4 =1$. Here, the holonomy $\alpha$ is not unitary but, the f\/irst author proved in~\cite{nonuni} that one can def\/ine a Floer cohomology with non-unitary line bundles.
As in the picture, we can list up all strips which bound $L_u$ and $L$. Then, the similar technique to the one given in Section~\ref{chanleungcp1} will give the corresponding holomorphic functions.
Note below that there are two more discs $D_5$ and $D_6$ which are not easily visible in  Fig.~\ref{dsteardrop}. We will be able to f\/ind these in the development picture (see Fig.~\ref{develop}).

\begin{enumerate}\itemsep=0pt
\item[(1)]
strips from $a$ to $b$:

\centerline{\includegraphics[height=0.5in]{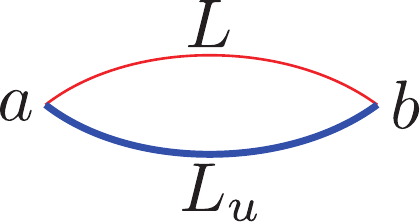}}
\vspace{-1mm}

\begin{enumerate}\itemsep=0pt
\item[$(i)$] The disc $D_1$ in the pictures leads to the term $-\frac{ z }{\alpha q}$.
\item[$(ii)$] In the limit, $D_2$ degenerates so that $D_2$ corresponds to $1$.
\end{enumerate}

\item[(2)] strips from $b$ to $a$

\centerline{\includegraphics[height=0.5in]{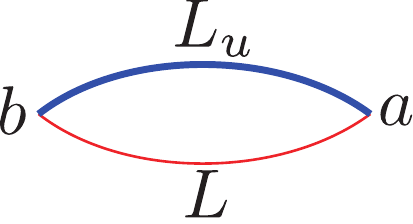}}

\begin{enumerate}\itemsep=0pt
\item[$(iii)$] $D_3$ gives rise to the term $\frac{q^4}{z} $.
\item[$(iv)$]  Since we only consider smooth disc in this subsection, $D_4$ in the picture indeed wraps around the singular point three times to give $- q^3\alpha^3 $. (We will consider nontrivial orbi-discs in the next subsection.)
\item[$(v)$] $D_5$ can be seen in the development f\/igure below which is clearly smooth since it covers the cone at the singular point three times. It leads to $- z q^2 \alpha^2 $.
\item[(vi)] Likewise, $D_6$ corresponds to the term~$- z^2 q \alpha $.
\end{enumerate}
\end{enumerate}

In conclusion,
\begin{gather}\label{MFtear}
\left( 1 - \dfrac{ z}{\alpha q}  \right) \left( \dfrac{q^4}{z} - q^3 \alpha^3 -  z q^2 \alpha^2 - z^2 q \alpha \right) = z^3 + \dfrac{q^4}{z} - \left(q^3 \alpha^3 + \dfrac{q^3}{\alpha} \right).
\end{gather}
By def\/inition of $\alpha$,  $q^3 \alpha^3 + \frac{q^3}{\alpha}$ is a critical value of $z^3 + \frac{q^4}{z}$.
\begin{prop}
The Lagrangian torus fiber for $u=0$, corresponds to the matrix factoriza\-tion~\eqref{MFtear} of $W= z^3 + \frac{q^4}{z}$ with critical value $\lambda =q^3 \alpha^3 + \frac{q^3}{\alpha}$.
\end{prop}

\subsection{The case of bulk deformed orbi-potential}\label{turnonbulk}
Now, we turn on bulk deformation by twisted sectors. Namely, for $\nu =[1] \in \Z/3$,
we can consider $X_\nu$ which is an isolated point, whose fundamental class is $1_\nu \in
H^0(X_\nu)$. Then we take $\frak{b} =c 1_\nu$, with $c \in \Lambda_{+}$.
Here,
\[
\Lambda_+ = \left\{ \left.\sum a_i T^{\lambda_i} \in \Lambda \, \right| \lambda_i >0 \right\}.
\]

Since there is an insertion from twisted sectors, now we also include orbifold holomorphic discs (orbi-discs for short). (We refer readers to~\cite{CP} for details of the following constructions. See~\cite{FOOO3}
for bulk deformations in the case of toric manifolds.)

We f\/irst f\/ind the bulk-deformed mirror.
As we have chosen $\frak{b} =c 1_\nu$, we need to consider orbifold holomorphic discs with
several orbifold interior marked points with $\Z/3$ singularity, mapping to $X_\nu$ where each generator of the local group at orbifold marked point is mapped to~$\nu$.
By simple degree consideration, the orbi-disc with only one orbifold interior marked point contributes to the potential, and such holomorphic orbi-disc is classif\/ied in~\cite{CP}. In this case, there is a unique
holomorphic orbi-disc $D_4'$, which covers the cone once.
The additional information from the orbi-disc(which has area $u+1/3$) can be described as follows.
\[
e^{3x}T^{3u+1} + e^{-x}T^{1-u} + ce^{x}T^{u+1/3}
\]
or we can write $z = e^xT^{u + 1/3}$ and $T =q$, which gives
\[
W^{\frak b}= z^3 + \frac{q^{4/3}}{z} + cz.
\]

We remark that the bulk deformation which makes the f\/iber $L_u$ for any $u < 1/3$ to have a~non-trivial Floer cohomology
is $c = T^{2/3 -2u} - 3T^{2u + 2/3}$ or
$c = q^{2/3 - 2u} -3q^{2u + 2/3}$. (If $u \geq 1/3$, then $c \notin \Lambda_+$. In
fact, f\/ibers for $1> u > 1/3$, can be displaced from itself by using the open set obtained by
removing the cone point.)

The critical point equation  for $W^{\frak b}$ is
\[
3z^2  - \frac{q^{4/3}}{z^2} + c =0,
\]
whose solution is denoted as $q^{1-u}\alpha$.
Then the critical value of the bulk deformed potential is
\[
q^{3-3u}\alpha^3 + q^{1/3 + u}/\alpha + cq^{1-u} \alpha.
\]

Now, we look at the corresponding matrix factorization.

We repeat \eqref{MFtear}, with the additional orbifold holomorphic strip contribution (underlined term below),
which is $D_4'$ in the development picture (Fig.~\ref{develop}), to obtain the following
matrix factorization:
\begin{gather}
 \left( 1 - \dfrac{ z}{q^{1-u}\alpha}  \right) \left( \dfrac{q^{4/3}}{z} - q^{3-3u}\alpha^3 -  z q^{2-2u}\alpha^2 - z^2 q^{1-u}\alpha - \underline{ c q^{1-u}\alpha} \right) \nonumber\\
\qquad{}= z^3 + \dfrac{q^{4/3}}{z} + cz - \left(q^{3-3u}\alpha^3 + q^{1/3 + u}/\alpha + cq^{1-u} \alpha \right).\label{MFtearorb1}
\end{gather}
\begin{prop}
For the teardrop orbifold $X$, with the bulk deformation $\frak{b} =c 1_\nu \in H^0(X_\nu)$, the Lagrangian fiber
$L_u$ for $u < 1/3$ corresponds to the matrix factorization~\eqref{MFtearorb1} of the
bulk deformed potential $W^{\frak b}$.
\end{prop}

\section{Weighted projective lines}\label{SecWPR}

Finally, we study the case of weighted $\CC P^1$s with general weights at ends. Let $X$ be a weighted projective line with $\ZZ/m\ZZ$-singularity on the left and $\ZZ / n\ZZ$ on the right, i.e.~$X$ is obtained by dividing $\CC^2 \setminus \{ (0,0)\}$ by the following action of $\CC^\ast$ ($m$, $n$ are assumed to be relatively prime):
\[
\rho \in \CC^\ast : \ (z,w) \mapsto (\rho^n z, \rho^m w).
\]

Recall that $\left[-\frac{1}{m},\frac{1}{n} \right]$ is the moment polytope of $X$.
(Similarly, one can also work on the case of toric orbifold of dimension one, which
corresponds to the interval as a polytope with integer labels $m$, $n$ at each end point.
In this setting, $m$, $n$ do not need to be relatively prime. But we leave the details of this general case to the
interested reader.)

 Thus, in this case, smooth holomorphic discs of Maslov index 2 can be described as follows:
\[
e^{mx} T^{1+mu} + e^{-nx} T^{1-nu}.
\]
As before, we make a substitution for this equation by $z= e^{x} T^{\frac{1}{m}+u}$ and $q = T^{1/m}$, which is coherent to the computation for the teardrop. Therefore, the LG potential $W$ is given as follows:
\[
W=z^m+\dfrac{q^{m+n}}{z^n}.
\]
In terms of the {\it smooth} Floer theory, the f\/iber at $u=0$ has a non-vanishing Floer homology.
This is because two smooth holomorphic discs (left, and right) has the same area at $u=0$.

\begin{figure}
\centering
\includegraphics[height=2.4in]{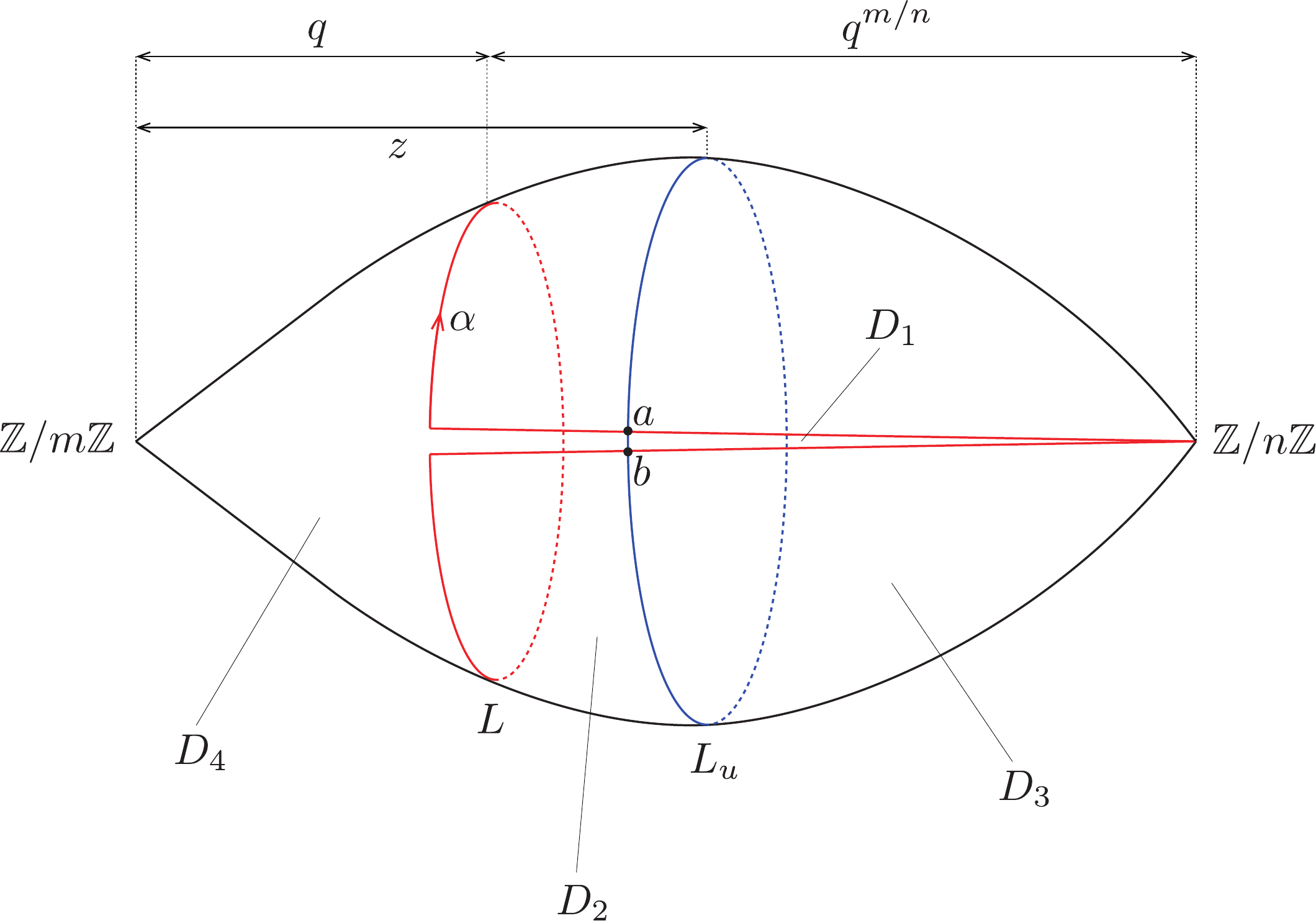}
\caption{Disk splittings in weighted projective lines.}
\label{fig:weightproj}
\end{figure}

As before, we consider the deformation $L$ of the central f\/iber $L_0$ and consider intersection with a general f\/iber $L_u$ at $u \in \left( 0,\frac{1}{n} \right) \subset \left[-\frac{1}{m},\frac{1}{n} \right].$
The total area of $X$ corresponds to{\samepage
\[
T^{\frac{1}{n} - \left(-\frac{1}{m}\right)} = q^{\frac{m+n}{n}},
\]
and $L$ splits the total area into $q$ and $q^{m/n}$  (see Fig.~\ref{fig:weightproj}).}

Let $\alpha$ be a (non-unitary) holonomy of $L$, which is given by one of the solutions of the equation
\begin{gather}\label{weightedholo}
m z^{m-1}-\dfrac{n}{z^{n+1}}=0.
\end{gather}

Then we count index 2 holomorphic strips as we did above, counting visible ones from Fig.~\ref{fig:weightproj} and those from development f\/igures which are obtained by letting singularities at both ends being~$\infty$. (Strips can cover the region of $D_3$ and $D_4$ several times.)
We remark that only~$D_1$ and~$D_2$ are strips from~$a$ to~$b$, and any other strips such as~$D_3$, $D_4$ and those given by development f\/igures are strips from~$b$ to~$a$.

After counting all such strips, we have the factorization as
\begin{gather}
 \left( 1-\dfrac{z}{\alpha q} \right) \left( \displaystyle\sum_{k=0}^{n} \frac{q^{ \frac{m}{n}  k}}{\alpha^{k}} \left(\frac{q^{\frac{m+n}{n}}}{z} \right)^{n-k} -\sum_{k=1}^{m}\alpha^k q^{k}z^{m-k} \right) \nonumber \\
\qquad{}= z^m+\dfrac{q^{m+n}}{z^n}-\left( \alpha^m q^{m}+\dfrac{q^{m}}{\alpha^n} \right).\label{MFweighted}
\end{gather}
(The f\/irst factor in the left-hand side of \eqref{MFweighted} counts the strips from $a$ to $b$.)
One can easily see that $\alpha^m q^{m}+\frac{q^{m}}{\alpha^n}$ is a critical value of $W=z^m+\frac{q^{(m+n)}}{z^n}$, comparing with~\eqref{weightedholo}.  Note that if $m=3$ and $n=1$, then the result coincides with the one we have obtained in the previous section.

\subsection*{Acknowledgements}

First author thank  Naichung Conan Leung, Siu Cheong Lau, Kwok Wai Chan, Yong-Geun Oh for
helpful discussions and for sharing their ideas. The work of C.-H.~Cho was supported by the National Research Foundation of Korea Grant funded by the Korean Government MEST, Basic Research Promotion Fund (NRF-2011-013-C0004).

\pdfbookmark[1]{References}{ref}
\LastPageEnding


\begin{thebibliography}{99}
\footnotesize\itemsep=0pt

\bibitem{Al}
Alston G., Floer cohomology of real Lagrangians in the Fermat quintic
  threefold, \href{http://arxiv.org/abs/1010.4073}{arXiv:1010.4073}.

\bibitem{ADD}
Ashok S.K., Dell'Aquila E., Diaconescu D.E., Fractional branes in
  {L}andau--{G}inzburg orbifolds, \textit{Adv. Theor. Math. Phys.} \textbf{8}
  (2004), 461--513, \href{http://arxiv.org/abs/hep-th/0401135}{hep-th/0401135}.

\bibitem{Au}
Auroux D., Mirror symmetry and {$T$}-duality in the complement of an
  anticanonical divisor, \textit{J.~G\"okova Geom. Topol. GGT} \textbf{1}
  (2007), 51--91, \href{http://arxiv.org/abs/0706.3207}{arXiv:0706.3207}.

\bibitem{CL}
Chan K., Leung N.C., Matrix factorizations from SYZ transformations, in
  Advances in Geometric Analysis, \textit{Adv. Lect. Math.}, Vol.~21,
  International Press, Somerville, MA, 2011, 203--224, \href{http://arxiv.org/abs/1006.3832}{arXiv:1006.3832}.

\bibitem{CL2}
Chan K., Leung N.C., Mirror symmetry for toric {F}ano manifolds via {SYZ}
  transformations, \href{http://dx.doi.org/10.1016/j.aim.2009.09.009}{\textit{Adv. Math.}} \textbf{223} (2010), 797--839,
  \href{http://arxiv.org/abs/0801.2830}{arXiv:0801.2830}.


\bibitem{C2}
Cho C.-H., Constant triangles in Fukaya category, in preparation.

\bibitem{C1}
Cho C.-H., Holomorphic discs, spin structures, and {F}loer cohomology of the
  {C}lif\/ford torus, \href{http://dx.doi.org/10.1155/S1073792804132716}{\textit{Int. Math. Res. Not.}} \textbf{2004} (2004), no.~35,
  1803--1843, \href{http://arxiv.org/abs/math.SG/0308224}{math.SG/0308224}.

\bibitem{nonuni}
Cho C.-H., Non-displaceable {L}agrangian submanifolds and {F}loer cohomology
  with non-unitary line bundle, \href{http://dx.doi.org/10.1016/j.geomphys.2008.06.003}{\textit{J.~Geom. Phys.}} \textbf{58} (2008),
  1465--1476, \href{http://arxiv.org/abs/0710.5454}{arXiv:0710.5454}.

\bibitem{CO}
Cho C.-H., Oh Y.-G., Floer cohomology and disc instantons of {L}agrangian torus
  f\/ibers in {F}ano toric manifolds, \textit{Asian~J. Math.} \textbf{10} (2006),
  773--814, \href{http://arxiv.org/abs/math.SG/0308225}{math.SG/0308225}.

\bibitem{CP}
Cho C.-H., Poddar M., Holomorphic orbi-discs and Lagrangian Floer cohomology for
  toric orbifolds, \href{http://arxiv.org/abs/1206.3994}{arXiv:1206.3994}.



\bibitem{FOOO2}
Fukaya K., Oh Y.-G., Ohta H., Ono K., Lagrangian {F}loer theory on compact toric
  manifolds.~{I}, \href{http://dx.doi.org/10.1215/00127094-2009-062}{\textit{Duke Math.~J.}} \textbf{151} (2010), 23--174,
  \href{http://arxiv.org/abs/0802.1703}{arXiv:0802.1703}.

\bibitem{FOOO3}
Fukaya K., Oh Y.-G., Ohta H., Ono K., Lagrangian {F}loer theory on compact toric
  manifolds. {II}.~Bulk deformations, \href{http://dx.doi.org/10.1007/s00029-011-0057-z}{\textit{Selecta Math.~(N.S.)}} \textbf{17}
  (2011), 609--711, \href{http://arxiv.org/abs/0810.5654}{arXiv:0810.5654}.

  \bibitem{FOOO}
Fukaya K., Oh Y.-G., Ohta H., Ono K., Lagrangian intersection {F}loer theory:
  anomaly and obstruction, \textit{AMS/IP Studies in Advanced Mathematics},
  Vol.~46, Amer. Math. Soc., Providence, RI, 2009.

\bibitem{Gross}
Gross M., The {S}trominger--{Y}au--{Z}aslow conjecture: from torus f\/ibrations
  to degenerations, in Algebraic Geometry~-- {S}eattle 2005, \textit{Proc.
  Sympos. Pure Math.}, Vol.~80, Amer. Math. Soc., Providence, RI, 2009,
  {P}art~1, 149--192.

\bibitem{KL}
Kapustin A., Li Y., D-branes in {L}andau--{G}inzburg models and algebraic
  geometry, \href{http://dx.doi.org/10.1088/1126-6708/2003/12/005}{\textit{J.~High Energy Phys.}} \textbf{2003} (2003), no.~12, 005, 44~pages, \href{http://arxiv.org/abs/hep-th/0210296}{hep-th/0210296}.

\bibitem{KO}
Kwon D., Oh Y.-G., Structure of the image of (pseudo)-holomorphic discs with
  totally real boundary condition, \textit{Comm. Anal. Geom.} \textbf{8}
  (2000), 31--82.

\bibitem{Oh}
Oh Y.-G., Floer cohomology of {L}agrangian intersections and pseudo-holomorphic
  disks.~{I}, \href{http://dx.doi.org/10.1002/cpa.3160460702}{\textit{Comm. Pure Appl. Math.}} \textbf{46} (1993), 949--993.

\bibitem{Oh2}
Oh Y.-G., Floer cohomology, spectral sequences, and the {M}aslov class of
  {L}agrangian embeddings, \href{http://dx.doi.org/10.1155/S1073792896000219}{\textit{Int. Math. Res. Not.}} \textbf{1996} (1996),
no.~7, 305--346.

\bibitem{Or}
Orlov D.O., Triangulated categories of singularities and {D}-branes in
  {L}andau--{G}inzburg models, \textit{Tr. Mat. Inst. Steklova} \textbf{246}
  (2004), 240--262.

\bibitem{S}
Seidel P., Fukaya categories and {P}icard--{L}efschetz theory, \href{http://dx.doi.org/10.4171/063}{\textit{Z\"urich
  Lectures in Advanced Mathematics}}, European Mathematical Society (EMS),
  Z\"urich, 2008.

\bibitem{SYZ}
Strominger A., Yau S.T., Zaslow E., Mirror symmetry is $T$-duality, in Winter
  {S}chool on {M}irror {S}ymmetry, {V}ector {B}undles and {L}agrangian
  {S}ubmanifolds ({C}ambridge, 1999), \textit{AMS/IP Studies in Advanced
  Mathema\-tics}, Vol.~23, Editors C.~Vafa, S.T.~Yau, Amer. Math. Soc.,
  Providence, RI, 2001, 275--295.

\end{thebibliography}
\end{document}